\documentclass[12pt]{article}
\usepackage{amscd, amssymb, amsmath}
\begin{document}
\renewcommand{\thesubsection}{\arabic{subsection}}
\newenvironment{proof}{{\bf Proof}:}{\vskip 5mm }
\newtheorem{proposition}{Proposition}[subsection]
\newtheorem{lemma}[proposition]{Lemma}
\newtheorem{definition}[proposition]{Definition}
\newtheorem{theorem}[proposition]{Theorem}
\newtheorem{cor}[proposition]{Corollary}
\newtheorem{conjecture}[proposition]{Conjecture}
\newtheorem{pretheorem}[proposition]{Pretheorem}
\newtheorem{hypothesis}[proposition]{Hypothesis}
\newtheorem{example}[proposition]{Example}
\newtheorem{remark}[proposition]{Remark}
\newcommand{\comment}[1]{}
\newcommand{\mcal}{\mathcal}
\newcommand{\sr}{\rightarrow}
\newcommand{\lr}{\longrightarrow}
\newcommand{\zz}{{\bf Z}}
\newcommand{\zq}{{\bf Z}_{qfh}}
\newcommand{\nn}{{\bf N}}
\newcommand{\qq}{{\bf Q}}
\newcommand{\nq}{{\bf N}_{qfh}}
\newcommand{\oo}{\otimes}
\newcommand{\uu}{\underline}
\newcommand{\ih}{\uu{Hom}}
\newcommand{\af}{{\bf A}^1}
\newcommand{\ro}{\mu}
\newcommand{\llabel}[1]{\label{#1}}
\newcommand{\BB}{\bullet}
\newcommand{\wt}{\widetilde}

\begin{center}
{\Large\bf On 2-torsion in motivic cohomology}\\
\vskip 4mm {\large\bf Vladimir Voevodsky}\footnote{The work on this paper
was supported by NSF grants DMS-96-25658, DMS-97-29992 and DMS-99-01219,
Abrose Monell Foundation, Clay Mathematics Institute, Sloan Research
Fellowship and Veblen Fund}$^,$\footnote{School of Mathematics, Institute for Advanced
Study, Princeton NJ, USA. e-mail: vladimir@ias.edu}\\
\vskip 3mm
{June 2001}
\end{center}
\vskip 25mm
\tableofcontents
\subsection{Introduction}
In the early 1980-ies Alexander Beilinson and Stephen Lichtenbaum (see
\cite{Licht2}, \cite[5.10.D]{pairing}) introduced the idea of motivic
complexes $\zz(n)$ and formulated a set of conjectures describing
their properties.  Among these conjectures was the following one which
they called the generalized Hilbert 90 property.
\begin{conjecture}
\llabel{realmain} For any field $k$ and any $n\ge 0$ one has 
$${\bf
H}^{n+1}_{et}(Spec(k),\zz(n))=0$$
\end{conjecture}
As we know today it implies an amazing number of other conjectures
about Galois cohomology, motivic cohomology and algebraic
K-theory. The goal of this paper is to prove the 2-local version of
Conjecture \ref{realmain}.
\begin{theorem}
\llabel{realmainth}
For any field $k$ and any $n\ge 0$ one has 
$${\bf H}^{n+1}_{et}(Spec(k),\zz_{(2)}(n))=0$$
\end{theorem}
As a corollary of Theorem \ref{realmainth} we get the following result.
\begin{theorem}
\llabel{milnorcon}
Let $k$ be a field of characteristic not equal to $2$. Then the norm
residue homomorphism
$$K_*^M(k)/2\sr H^*(k,\zz/2)$$
is an isomorphism.
\end{theorem}
The statement of Theorem \ref{milnorcon} is known as the Milnor
Conjecture or the $\zz/2$-coefficients case of the Bloch-Kato
conjecture. The general ``Bloch-Kato Conjecture'' can be described as
follows. Let $k$ be a field and $l$ a prime number different from the
characteristic of $k$. Fix a separable closure $k_{sep}$ of $k$ and
let $\mu_l$ denote the group of $l$-th roots of unity in $k_{sep}$. One
may consider $\mu_l$ as a $Gal(k_{sep}/k)$-module. By definition of
$\mu_l$ one has a short exact sequence
$$1\lr \mu_l\lr k_{sep}^*\stackrel{z^l}{\lr} k_{sep}^*\lr 1$$
which is called the Kummer sequence. The boundary map in the
associated long exact sequence of Galois cohomology is a homomorphism
\begin{equation}
\llabel{start}
k^*\sr H^1(k,\mu_l)
\end{equation}
In \cite{BT}, Bass and Tate proved that for $a\in k^*-\{1\}$ the
cohomology class $(a)\wedge (1-a)$ lying in $H^2(k,\mu_l^{\oo 2})$ is
zero i.e. that the homomorphism (\ref{start}) extends to a
homomorphism of rings
\begin{equation}
\llabel{start2}
T(k^*)/I\sr H^*(k,\mu_l^{\oo *})
\end{equation}
where $T(k^*)$ is the tensor algebra of the abelian group $k^*$ and
$I$ the ideal generated by elements of the form $a\oo b$ for $a,b\in
k^*$ such that $a+b=1$. The graded components of the quotient
$T(k^*)/I$ are known as the Milnor K-groups of $k$ and the
homomorphism (\ref{start2}) is usually written as
\begin{equation}
\llabel{start3}
K_*^M(k)\sr H^*(k,\mu_l^{\oo *})
\end{equation}
Let $K_*^M(k)/l$ be the quotient of $K_*^M(k)$
by the ideal of elements divisible by $l$. It is clear that
(\ref{start3}) factors through a map 
\begin{equation}
\llabel{start4}
K_*^M(k)/l\sr H^*(k,\mu_l^{\oo *})
\end{equation}
which is called the {\em norm residue homomorphism}. Theorem
\ref{milnorcon} proves the following conjecture for $l=2$.
\begin{conjecture}
\llabel{BK}
The map (\ref{start4}) is an isomorphism for any field $k$ of
characteristic $\ne l$.
\end{conjecture}
This conjecture has a long and convoluted history.  The map
(\ref{start4}) is clearly an isomorphism in degree zero. In degree
one, (\ref{start4}) is a monomorphism and its cokernel is the group of
$l$-torsion elements in $H^1(k,k_{sep}^*)$ which is known to be zero
as a corollary of the classical Hilbert 90 Theorem.

In degree 2 the homomorphism (\ref{start4}) has an interpretation in
terms of central simple algebras. The question about central simple
algebras which is equaivalent to surjectivity of (\ref{start4}) in
degree 2 seems to be very old. The question of injectivity of
(\ref{start4}) in degree 2 was explicitly stated by John Milnor in
\cite[Remark on p. 147]{Milnor2}.
 
In 1981 Alexander Merkurjev published a paper \cite{Merkurjev} where
he proved that (\ref{start4}) is an isomorphism in degree 2 for $l=2$
and any field $k$ (such that $char(k)\ne 2$). This paper is the
starting point of all the further work on bijectivity of
(\ref{start4}). In 1982 Merkurjev together with Andrei Suslin proved
that (\ref{start4}) is an isomorphism in degree $2$ for all $l$ (see
\cite{MS1}).

In degree 3 and $l=2$ the bijectivity of (\ref{start4}) was proved by
Merkurjev and Suslin in \cite{MS2} and independently by Markus Rost in
\cite{Rost3}. 

In \cite{Milnor} Milnor considered the homomorphism (\ref{start4}) in
all degrees for $l=2$ as a part of his investigation of the relations
between $K_*^M(k)/2$ and quadratic forms over $k$. He mentioned that
he does not know of any fields $k$ for which (\ref{start4}) fails to
be an isomorphism and gave several examples of classes of fields $k$
for which (\ref{start4}) is an isomorphism in all degrees. His
examples extend a computation made by Bass and Tate in the same paper
where they introduce (\ref{start}).  This paper of Milnor is the
reason why Conjecture \ref{BK} for $l=2$ is called the Milnor
Conjecture.

The name ``Bloch-Kato conjecture'' comes from a 1986 paper by Spencer
Bloch and Kazuya Kato \cite{BK} where they mention (on p. 118) that
``one conjectures [the homomorphism (\ref{start4})] to be an
isomorphism quite generally''. This name is certainly incorrect
because Conjecture \ref{BK} appeared quite explicitly in the work of
Beilinson \cite{pairing} and especially Lichtenbaum \cite{Licht2}
which was published much earlier.

The list of hypothetical properties of motivic complexes of Beilinson
and Lichtenbaum implied that
\begin{equation}
\llabel{equ1}
{\bf H}^{n}_{Zar}(Spec(k),\zz/l(n))=K^n_M(k)/l
\end{equation}
and if $l\ne char(k)$ then 
\begin{equation}
\llabel{equ2}
{\bf H}^{n}_{et}(Spec(k),\zz/l(n))=H^n(k,\mu_l^{\oo n})
\end{equation}
The norm residue homomorphism in this language is just the canonical
map from the motivic cohomology in the Zariski topology to the motivic
cohomology in the etale topology. The conjectures made by Beilinson
and Lichtenbaum implied the following:
\begin{conjecture}
\llabel{blmod} Let $X$ be a smooth variety over a field $k$. Then the
map
$${\bf H}^{p}_{Zar}(X,\zz(q))\sr {\bf H}^{p}_{et}(X,\zz(q))$$
from the motivic cohomology of $X$ in the Zariski topology to the
motivic cohomology of $X$ in the etale topology is an isomorphism for
$p\le q+1$.
\end{conjecture}
In view of (\ref{equ1}) and (\ref{equ2}) Conjecture \ref{blmod}
implies Conjecture \ref{BK}.  Note also that since for any field $k$
one has $H^{n+1}_{Zar}(Spec(k),\zz(n))=0$, Conjecture \ref{realmain}
is a particular case of Conjecture \ref{blmod}.

\vskip 3mm
\noindent
The relation between Conjecture \ref{blmod} and Conjecture \ref{BK}
was further clarified in \cite{SusVoe3} and \cite{GL} where it was
shown that the $l$-local version of Conjecture \ref{blmod} for fields
$k$ such that $char(k)\ne l$ is in fact equivalent to Conjecture
\ref{BK} and moreover that it is sufficient to show only the
surjectivity of (\ref{start4}), because the injectivity formally
follows. In Theorem \ref{sv} we show using results of \cite{SusVoe3}
and \cite{GL} that the $l$-local version of Conjecture \ref{realmain}
implies the $l$-local version of Conjecture \ref{blmod} for
$char(k)\ne l$.

\vskip 3mm
\noindent
In \cite{GL2} Thomas Geisser and Marc Levine proved the p-local
version of Conjecture \ref{blmod} for schemes over fields of
characterstic $p$. They use the version of motivic cohomology based on
the higher Chow groups which is now known to be equivalent to the
version used here by \cite{Suslin3new}, \cite{FS} and
\cite{comparison}. Together with our Theorem \ref{sv} this result
implies in particular that Conjecture \ref{blmod} is a corollary of
Conjecture \ref{realmain} for all $k$.

\vskip 3mm
\noindent
The first version of the proof of Theorem \ref{realmainth} appeared in
\cite{MC0}. It was based on the idea that there should exist algebraic
analogs of the higher Morava K-theories and that the $m$-th algebraic
Morava K-theory can be used for the proof of Conjecture \ref{BK} for
$l=2$ and $n=m+2$ in the same way as the usual algebraic K-theory is
used in Merkurjev-Suslin proof of Theorem \ref{milnorcon} in degree
$3$ in \cite{MS2}. This approach was recently validated by Simone
Borghesi \cite{Simone} who showed how to construct algebraic Morava
K-theories (at least in characteristic zero).

The second version of the proof appeared in \cite{MC}. Instead of
algebraic Morava K-theories it used small pieces of these theories
which are easy to construct as soon as one knows some facts about the
cohomological operations in motivic cohomology and their
interpretation in terms of the motivic stable homotopy category.

The main difference between the present paper and \cite{MC} is in the
proof of Theorem \ref{newth2} (\cite[Theorem 3.25]{MC}). The approach
used now was outlined in \cite[Remark on p.39]{MC}. It is based on the
connection between cohomological operations and characteristic classes
and circumvents several technical ingredients of the older proof. The
most important simplification is due to the fact that we can now
completely avoid the motivic {\em stable} homotopy category and the
topological realization functor.

Another difference between this paper and \cite{MC} is that we can now
prove all the intermediate results for fields of any
characteristic. Several developments made this possible. The new proof
of the suspension theorem for the motivic cohomology \cite[Theorem
2.4]{Red} based on the comparison between the motivic cohomology and
the higher Chow groups established in \cite{Suslin3new}, \cite{FS} and
\cite{comparison} does not use resolution of singularities. The same
comparison together with the new proof of the main result of
\cite{SusVoe3} by Thomas Geisser and Marc Levine in \cite{GL} allows
one to drop the resolution of singularities assumption in the proof of
Theorem \ref{sv}. Finally, the new approach to the proof of Theorem
\ref{newth2} does not require the topological realization functor
which only exists in characteristic zero.

\vskip 2mm
\noindent
The paper is organized as follows. In Sections \ref{sec1}-\ref{sec4}
we prove several results about motivic cohomology which are used in
the proof of our main theorem but which are not directly related to
the Belinson-Lichtenbaum conjecturs. In Section \ref{sec2} this is
Corollary \ref{exitcor}, in Section \ref{sec3} Theorems \ref{Pmain}
and \ref{rostvan} and in Section \ref{sec4} Theorem \ref{l4}. The
proof of Corollary \ref{exitcor} uses Theorem \ref{newth8} of Section
\ref{sec1}. There are no other connections of these four sections to
each other or to the remaining sections of the paper.

In Section \ref{sec5} we show that Conjecture \ref{blmod} is a
corollary of Conjecture \ref{realmain} and that Conjecture \ref{BK} is
a corollary of Conjecture \ref{blmod}. This section is independent of
the previous four sections.

In Section \ref{sec6} we prove Theorem \ref{realmainth}. This section
also contains some corollaries of the main theorem. Using similar
techniques one can also prove the Milnor Conjecture which asserts that
the Milnor ring modulo 2 is isomorphic to the graded Witt ring of
quadratic forms. For the proof of this result together with a more
detailed computations of motivic cohomology groups of norm quadrics
see \cite{OVV}.

Two appendixes contain the material which is used throughout the paper
and which I could not find good references for.

\vskip 3mm
\noindent
All through the paper we use the Nisnevich topology \cite{Nisnevich}
instead of the Zariski one. Since all the complexes of sheaves
considered in this paper have transfers and homotopy invariant
cohomology sheaves \cite[Theorem 5.7, p.128]{H2new} implies that one
can replace Nisnevich hypercohomology by Zariski ones everywhere in
the paper without changing the answers.

\vskip 3mm
\noindent
I am glad to be able to use this opportunity to thank all the people
who answered a great number of my questions during my work on the
Beilinson-Lichtenbaum conjectures. First of all I want to thank Andrei
Suslin who taught me the techniques used in
\cite{MS1},\cite{MS2}. Quite a few of the ideas of the first part of
the paper are due to numerous conversations with him. Bob Thomason
made a lot of comments on the preprint \cite{MC0} and in particular
explained to me why algebraic K-theory with $\zz/2$-coefficients has
no multiplicative structure, which helped to eliminate the assumption
$\sqrt{-1}\in k^*$ in Theorem \ref{maintheorem}.  Jack Morava and Mike
Hopkins answered a lot of my (mostly meaningless) topological
questions and I am in debt to them for not being afraid of things like
the Steenrod algebra anymore. The same applies to Markus Rost and
Alexander Vishik for not being afraid anymore of the theory of
quadratic forms. Dmitri Orlov guessed the form of the distinguished
triangle in Theorem \ref{Pmain} which was a crucial step to the
understanding of the structure of motives of Pfister quadrics. I would
also like to thank Fabien Morel, Chuck Weibel, Bruno Kahn and Rick
Jardine for a number of discussions which helped me to finish this
work. Finally, I would like to thank Eric Friedlander who introduced
me to Anderei Suslin and helped me in many ways during the years when I
was working on the Theorem \ref{milnorcon}.

Most of the mathematics of this paper was invented when I was a Junior
Fellow of the Harvard Society of Fellows and I wish to express my deep
gratitude to the society for providing a unique opportunity to work
for three years without having to think of things earthly. The first
complete version was written during my stay in the Max-Planck
Institute in Bonn. Further work was done when I was at the
Northwestern University and in its final form the paper was written
when I was a member of the Institute for Advanced Study in Princeton.
\subsection{The degree map}
\llabel{sec1}
\llabel{degmap}
In this section and Section \ref{sec2} we use the motivic homotopy
theory of algebraic varieties developed in \cite{MoVo}. For an
introduction to this theory see also \cite{talk} and
\cite{delnotes}. For a smooth projective variety $X$ of pure dimension
$d$ we consider the degree map $H^{2d,d}(X,\zz)\sr \zz$ (``evaluation
on the fundamental class''). Most of the section is occupied by the
proof of Theorem \ref{newth8} where we show that the degree map can be
described as the composition of the Thom isomorphism for an
appropriate vector bundle and a map defined by a morphism in the {\em
pointed motivic homotopy category} $H_{\BB}$. This theorem is a formal
corollary of Spanier-Whitehead duality in the motivic stable
homotopy theory but since the details of this duality are not worked
out yet we are forced to give a very direct but not a very conceptual
proof here.

\vskip 3mm
\noindent
Let $X$ be a smooth projective variety of pure dimension $d$. As we
know from \cite{Suslin3new}, \cite{FS} and \cite{comparison} there is
an isomorphism $H^{2q,q}(X,\zz)=CH^q(X)$ where $CH^q(X)$ is the group
of rational equivalence classes of cycles of codimension $q$ on
$X$. In particular, $H^{2d,d}(X,\zz)=CH^0(X)$ and there is the degree
map $deg:H^{2d,d}(X,\zz)\sr \zz$. 

For a closed embedding $i:Z\sr X$ of smooth varieties over $k$ with
the normal bundle $N$ denote by $\gamma_i$ the composition of the
projection $X_+\sr X/(X-Z)$ with the weak equivalence $X/(X-Z)\sr
Th_Z(N)$ constructed in \cite[Theorem 2.23]{MoVo}. We will use
the following properties of the degree map.
\begin{lemma}
\llabel{degprop1} Let $i:X\sr Y$ be a closed embedding of smooth
projective schemes over $k$ of pure dimensions $d_X$ and $d_Y$
respectively. Let $N$ be the normal bundle to $i$ and $t_N$ the
Thom class of $N$. Then for any $a\in H^{2d_X, d_X}(X,\zz)$ one has:
$$deg(\gamma_i^*(a\wedge t_N))=deg(a)$$
\end{lemma}
\begin{lemma}
\llabel{newl7} Let $i:X\sr {\bf P}^d$ be a smooth hypersurface of
degree $n$. Then $deg(i^*(e({\cal O}(1))^{d-1}))=n$.
\end{lemma}
%
%\begin{proof}
%???
%\end{proof}
%
\begin{lemma}
\llabel{newl12} Let $X$ be a smooth projective variety of pure
dimension $d$. Then the map
$$H^{2d,d}(X,\zz)\stackrel{deg}{\sr}\zz\sr \zz/n$$ 
is surjective if and only if $X$ has a rational point over an
extension of $k$ of degree prime to $n$.
\end{lemma}
%
%\begin{proof}
%???
%\end{proof}
%
\begin{remark}\rm
For a smooth projective variety $X$ of dimension $d$
and elements $a\in H^{*,>0}(Spec(k),\zz)$ and $x\in H^{*,*}(X)$
such that $ax\in H^{2d,d}(X)$ one clearly has $deg(ax)=0$. Lemma
\ref{degprop1} together with this property and the property that
$deg(1)=1$ uniquely characterize the degree maps $H^{2d,d}(X,\zz)\sr
\zz$ for projective $X$. This can be easily seen using \cite[Lemma
4.2]{Red}.
\end{remark}
For the proof of the Theorem \ref{newth8} we will need the following
construction. Let $i:Z\sr X$ be a closed embedding of smooth varieties
over $k$ with the normal bundle $N$. Let further $V$ be a vector
bundle on $X$ and $X\sr V$ the zero section. Consider the
embedding $Z\sr X\sr V$. Since the normal bundle to $Z$ in $V$ is
$N\oplus i^*(V)$ we have a weak equivalence $V/(V-Z)\sr Th_Z(N\oplus
i^*(V))$. This gives us a map
$$\gamma_{i,V}:Th_X(V)\cong V/(V-X)\sr V/(V-Z)\cong Th_Z(N\oplus
i^*(V)).$$
\begin{lemma}
\llabel{easynew}
Let $x$ be a class in $H^{*,*}(Z)$. Then one has:
$$\gamma_{i,V}^*(xt_{N\oplus i^*(V)})=\gamma_i^*(xt_N)t_V$$
where $t_{N\oplus i^*(V)}$, $t_V$ and $t_N$ are the Thom classes of
the corresponding vector bundles.
\end{lemma}
%
%\begin{proof}
%???
%\end{proof}
%
\begin{theorem}
\llabel{newth8} Let $X$ be a smooth projective variety of pure
dimension $d$ over a field $k$. Then there exists an integer $n$ and a
vector bundle $V$ on $X$ of dimension $n$ such that:
\begin{enumerate}
\item $V+T_X={\cal O}^{n+d}$ in $K_0(X)$
\item There exists a morphism in $H_{\bullet}$ of the form
$f_V:T^{n+d}\sr Th_X(V)$ such that the map $H^{2d,d}(X)\sr \zz$
defined by $f_V$ and the Thom isomorphisms coincides with the degree
map.
\end{enumerate}
\end{theorem}
\begin{proof}
In Proposition \ref{part1} below we show that the statement holds for
$X={\bf P}^d$ and a vector bundle $V_d$ on ${\bf P}^d$ of dimension
$n_d$. Let $i:X\sr {\bf P}^{m}$ be a closed embedding and $N$ the
normal bundle to $i$. Consider the composition
$$T^{n_m+m}\stackrel{f_{V_m}}{\sr}Th_{{\bf
P}^m}(V_m)\stackrel{\gamma_{i,V_m}}{\sr} Th_{X}(i^*(V_m)\oplus N)$$
If $\gamma_i:{\bf P}^m\sr Th_X(N)$ is the canonical map and $x$ is a class
in $H^{*,*}(X,\zz)$ then by Lemma \ref{easynew} one has
$$\gamma_{i,V_m}^*(xt_{i^*(V_m)\oplus N})=\gamma^*_i(xt_{N})t_{V_m}$$
Let $t$ be the tautological class in
$\wt{H}^{2(n_m+m),n_m+m}(T^{n_m+m},\zz)$, then 
$$f_{V_m}^*\gamma_{i,V_m}^*(xt_{i^*(V_m)\oplus
N})=f_{V_m}^*(\gamma_i^*(xt_{N})t_{V_m})=deg(\gamma_i^*(xt_{N}))t=$$
$$=deg(x)t$$
where the second equality holds by Proposition \ref{part1} and the
third one by Lemma \ref{degprop1}. Finally, in $K_0(X)$ we have
$T_X=i^*(T_{{\bf P}^m})-N$ and therefore,
$$T_X+N+i^*(V_m)={\cal O}^{n_m+m}.$$
\end{proof}
\begin{proposition}
\llabel{part1} There exists a vector bundle $V_d$ of dimension $n_d$
on ${\bf P}^d$ and a morphism $f_V:T^{d+n}\sr Th_{{\bf P}^d}(V)$ in
$H_{\BB}$ such that the homomorphism
$$H^{2d,d}({\bf P}^d,\zz)\sr \zz$$
defined by the Thom isomorphisms and $f_V$ coincides with the degree
map and $V+T_{{\bf P}^d}={\cal O}^{n_d+d}$ in $K_0({\bf P}^d)$. 
\end{proposition}
\begin{proof}
Let $W=\Omega\oplus(\Omega\oo T)$ where $T$ is the tangent bundle on
${\bf P}^d$ and $\Omega$ is its dual.
The dimension of $W$ is $n=d^2+d$. 
\begin{lemma}
\llabel{suma}
$W+T={\cal O}^{d^2+2d}$ 
\end{lemma}
\begin{proof}
We have to show that $T+\Omega+\Omega\oo T={\cal O}^{d^2+2d}$.
Consider the standard exact sequence
$$0\sr \Omega\sr {\cal O}(-1)^{d+1}\sr {\cal O}\sr 0$$
(see \cite[Th. 8.13]{Hartshorn}) and its dual
$$0\sr {\cal O}\sr {\cal O}(1)^{d+1}\sr T\sr 0.$$
The first sequence implies that $\Omega+{\cal O}=(d+1){\cal O}(-1)$,
hence
$$T+\Omega\oo T=T\oo({\cal O}+\Omega)=(d+1)T(-1).$$
Form the second sequence $T(-1)=(d+1){\cal O}-{\cal O}(-1)$, hence
$$\Omega+T+\Omega\oo T=\Omega+(d+1)T(-1)=\Omega+(d+1)^2{\cal
O}-(d+1){\cal O}(-1)=(d^2+2d){\cal O}.$$
\end{proof}
Consider the incidence hyperplane $H$ in ${\bf P}^d\times {\bf P}^d$
where the first projective space is thought of as the projective space
of a vector space ${\cal O}^{d+1}$ and the second one as the
projective space of $({\cal O}^{d+1})^*$. The complement 
$$\wt{\bf P}^d={\bf P}^d\times {\bf P}^d-H$$
considered as a scheme over ${\bf P}^d$ by means of a projection
$p:\wt{\bf P}^d\sr {\bf P}^d$ is an affine bundle over ${\bf P}^d$. On
the other hand, the Segre embedding 
$$i_{d,d}:{\bf P}^d\times {\bf P}^d\sr {\bf P}^{d^2+2d}$$
gives $H$ as the divisor at infinity for an appropriate choice of the
intersecting hyperplane $H_{\infty}$. Therefore, $\wt{\bf P}^d$ is an
affine variety. This construction is known as the Jouanolou trick (see
\cite{Jouanolou}).

Consider the fiber product:
\begin{equation}
\llabel{neweq}
\begin{CD}
\wt{\bf P}^d\times {\bf P}^d @>h>> {\bf P}^d\times {\bf P}^d\\
@VvVV @VVpV\\
\wt{\bf P}^d @>p>> {\bf P}^d
\end{CD}
\end{equation}
The open embedding $\wt{\bf P}^d\sr {\bf P}^d\times {\bf P}^d$ defines
a section $s$ of $v$. Let $N$ be the normal bundle to the Segre
embedding $i_{d,d}$ and let $E$ be the normal bundle to $s$.
\begin{lemma}
\llabel{secondl}
There is an integer $m$ and an isomorphism of vector bundles
\begin{equation}
\llabel{isonew}
E\oplus N\oplus {\cal O}^m\cong p^*(W)\oplus {\cal O}^m
\end{equation}
on $\wt{\bf
P}^n$.
\end{lemma}
\begin{proof}
Since $\wt{\bf P}^d$ is affine, two vector bundles give the same class
in $K_0$ if and only if they become isomorphic after the addition of
${\cal O}^m$ for some $m\ge 0$. Therefore, to prove the lemma it is
sufficient to show that one has $E+N+p^*(T_{{\bf P}^d})={\cal
O}^{d^2+2d}$ in $K_0(\wt{\bf P}^d)$. Let $p':\wt{\bf P}^d\sr {\bf
P}^d$ be the second of the two projections. One can easily see from
the definition of $E$ that $E=(p')^*(T_{{\bf P}^d})$. Hence
$E+p^*(T_{{\bf P}^d})$ is the restriction to $\wt{\bf P}^d$ of the
tangent bundle on ${\bf P}^d\times {\bf P}^d$. On the other hand the
short exact sequence which defines $N$ shows that over $\wt{\bf P}^d$
we have
$$T_{{\bf
P}^d\times {\bf P}^d}+N={\cal O}^{d^2+d}$$
which finishes the proof of the lemma. 
\end{proof}
\begin{lemma}
\llabel{firstl}
There exists a map $f:T^{d^2+2d}\sr Th_{\wt{\bf P}^d}(E\oplus N)$ such
that the homomorphism $H^{2d,d}({\bf P}^d)\sr \zz$ defined by $p^*$,
$f^*$ and the Thom isomorphisms coincides with the degree map.
\end{lemma}
\begin{proof}
Consider the pointed sheaf
$$F={\bf P}^{d^2+2d}/(H_{\infty}\cup ({\bf
P}^{d^2+2d}-i_{d,d}({\bf P}^d\times {\bf P}^d)))$$
There is an obvious map $\eta$ from $T^{d^2+2d}\cong {\bf
P}^{d^2+2d}/H_{\infty}$ to $F$. We claim that $F$ is isomorphic in
$H_{\bullet}$ to $Th_{\wt{\bf P}^n}(E\oplus N)$ and that the
composition of $\eta$ with this isomorphism satisfies the condition of
the lemma.

The same reasoning as in the proof of the homotopy purity theorem in
\cite[Theorem 2.23]{MoVo}, shows that in $H_{\BB}$ one has
$$F\cong Th_{{\bf P}^d\times{\bf P}^d}(N)/Th_{H}(N).$$
Consider again the pull-back square (\ref{neweq}). The map
$$Th_{\wt{\bf
P}^d\times {\bf P}^d}(h^*(N))/Th_{h^{-1}(H)}(h^*(N))\sr Th_{{\bf P}^d\times{\bf
P}^d}(N)/Th_{H}(N)$$ 
is clearly an $\af$-weak equivalence. On the other hand one verifies
easily that $h^{-1}(H)$ is contained in $\wt{\bf P}^d\times {\bf
P}^d-s(\wt{\bf P}^d)$ and that this embedding is an $\af$-weak
equivalence. Hence we have a weak equivalence from
$$Th_{\wt{\bf P}^d\times {\bf P}^d}(h^*(N))/Th_{h^{-1}(H)}(h^*(N))$$
to
$$Th_{\wt{\bf P}^d\times {\bf P}^d}(h^*(N))/Th_{\wt{\bf P}^d\times {\bf
P}^d-s(\wt{\bf P}^d)}(h^*(N))$$
The later quotient is isomorphic to
$h^*(N)/(h^*(N)-s(\wt{\bf P}^d))$ and since the normal
bundle to $s$ in $h^*{N}$ is $N\oplus E$ we conclude by
\cite[Theorem 2.23]{MoVo} that it is weakly equivalent to $Th_{\wt{\bf
P}^d}(E\oplus N)$. This finishes the construction of the map
$f:T^{d^2+2d}\sr Th_{\wt{\bf
P}^d}(E\oplus N)$. It remains to check that for $a\in H^{2d,d}({\bf
P}^d)$ we have $f^*(at_{E\oplus N})=deg(a)t$. 

Consider the following diagram
$$
\begin{CD}
{\bf P}^{d^2+2d}_+@>>> T^{d^2+2d}\\ 
@V\gamma_{i_{d,d}}VV @VVV\\ 
Th_{{\bf P}^d\times{\bf P}^d}(N) @>>> Th_{{\bf P}^d\times{\bf
P}^d}(N)/Th_H(N)\\
@A\tilde{h}AA @AAA\\ 
Th_{\wt{\bf P}^d\times{\bf P}^d}(N) @>>> Th_{\wt{\bf P}^d\times{\bf
P}^d}(N)/Th_{h^{-1}(H)}(N)\\
@V\gamma_{s,N}VV @VVV\\ 
Th_{\wt{\bf P}^d}(E\oplus N) @>=>> Th_{\wt{\bf P}^d}(E\oplus N)
\end{CD}
$$
where the arrows going up are weak equivalences. One verifies
immediately that it commutes. Denote the upper horizontal map of the
diagram by $\phi$. Since this map coincides with $\gamma_{i_0}$ where
$i_0:Spec(k)\sr {\bf P}^{d^2+2d}$ is the embedding of a point (see
\cite[Lemma 4.2]{Red}), Lemma \ref{degprop1} shows that
$deg(\phi^*(t))=1$ where $t$ is the tautological class in
$\wt{H}^{2(d^2+2d), d^2+2d}(T^{d^2+2d})$. Therefore, it it sufficient to
show that for a class $a\in H^{*,*}({\bf P}^d)$ we have
\begin{equation}
\llabel{neednew}
deg(\phi^*(f^*(p^*(a)t_{E\oplus
N})))=deg(a)
\end{equation}
Let $\tilde{h}$ denote the map
$$Th_{\wt{\bf P}^d\times{\bf P}^d}(N)\sr Th_{{\bf P}^d\times{\bf
P}^d}(N)$$
defined by $h$. Then
$$\phi^*(f^*(p^*(a)t_{E\oplus
N}))=\gamma_{i_{d,d}}^*(\wt{h}^{-1})^*\gamma^{*}_{s,N}(p^*(a)t_{E\oplus
N})=$$
$$=\gamma_{i_{d,d}}^*(\wt{h}^{-1})^*(\gamma_s^*(p^*(a)t_E)t_N)=\gamma_{i_{d,d}}^*((h^{-1})^*(\gamma_s^*(p^*(a)t_E))t_N)$$
The commutative diagram
$$
\begin{CD}
\wt{\bf P}^d\times {\bf P}^d @>\gamma_{s}>> Th_{\wt{\bf P}^d}(E)\\
@VhVV @VVV\\
{\bf P}^d\times {\bf P}^d @>\gamma_{\Delta}>> Th_{{\bf P}^d}(T_{{\bf
P}^d})
\end{CD}
$$
implies further that 
$$\gamma_s^*(p^*(a)t_E)=h^*\gamma_{\Delta}^*(at_{{T_{\bf P}^d}})$$
and therefore
$$\phi^*(f^*(p^*(a)t_{E\oplus
N}))=\gamma_{i_{d,d}}^*(\gamma_{\Delta}^*(at_{{T_{\bf P}^d}})t_N)$$
The equality (\ref{neednew}) follows now from Lemma \ref{degprop1}.
\end{proof}
We are ready now to finish the proof of Proposition \ref{part1}. We
take $V=W\oplus {\cal O}^m$ and $n=d^2+d+m$. The map $p:\wt{\bf
P}^d\sr {\bf P}^d$ defines a map of Thom spaces $Th_{\wt{\bf
P}^n}(p^*(V))\sr Th_{{\bf P}^n}(V)$. The isomorphism (\ref{isonew})
defines an isomorphism
$$\Sigma^m_TTh_{\wt{\bf P}^n}(E\oplus N)\cong Th_{\wt{\bf
P}^n}(p^*(V))$$
Composing these maps with the m-fold suspension of $f$ we get a map 
$$f_V:T^{d^2+2d+m}\sr Th_{{\bf P}^d}(V)$$
For $a\in H^{2d,d}({\bf P}^d)$ we have 
$$f_V^*(at_{V})=(\Sigma^m_Tf)^*(p^*(a)t_{E\oplus N\oplus{\cal
O}^m})=$$
$$=f^*(p^*(a)t_{E\oplus N})=deg(a)t$$
which shows that $f_V$ satisfies the condition of the proposition.
\end{proof}

\subsection{The motivic analog of Margolis
homology}\llabel{mainth}\llabel{sec2} 
In this section we introduce the motivic version of Margolis
homology. The definition of topological Margolis
homology\footnote{which appeared in \cite{Margolis} and which I
learned about from \cite{Rav1}.} is based on the fact that the
Steenrod algebra contains some very special elements $Q_i$ called
Milnor's primitives. These elements generate an exterior subalgebra in
the Steenrod algebra and in particular $Q_i^2=0$. Hence, one may
consider the cohomology of a space or a spectrum as a complex with the
differential given by $Q_i$. The homology of this complex are known in
topology as Margolis homology $\wt{MH}_{i}^*$. Spaces or spectra whose
Margolis homology vanish for $i\le n$ play an important role in the
proofs of the amazing recent results on the structure of the stable
homotopy category.  

Since we only know how to construct reduced power operations in the
motivic cohomology with $\zz/l$ coefficients for $l\ne char(k)$ 

\vskip 3mm
\noindent
{\em everywhere in this section $l$ is a prime not equal to the
characteristic of the base field.}

\vskip 3mm
\noindent
Recall that we defined in \cite[\S 13]{Red} operations $Q_i$ in the
motivic cohomology with $\zz/l$-coefficients of the form
$$Q_i:\wt{H}^{p,q}(-,\zz/l)\sr \wt{H}^{p+2l^i-1,q+l^i-1}(-,\zz/l)$$
We proved in \cite[Proposition 13.4]{Red} that operations $Q_i$ have
the property $Q_i^2=0$. For any ${\mcal X}$ the cohomology of the
complex
$$\wt{H}^{p-2l^i+1,q-l^i+1}({\mcal X},\zz/l)\sr
\wt{H}^{p,q}({\mcal X},\zz/l)\sr
\wt{H}^{p+2l^i-1,q+l^i-1}({\mcal X},\zz/l)$$
at the term $\wt{H}^{p,q}$ is called the $i$-th motivic Margolis
cohomology of ${\mcal X}$ of bidegree $(p,q)$ denoted by
$\wt{MH}_i^{p,q}({\mcal X})$. By \cite[Lemma 13.5]{Red} $Q_0$ is the
Bockstein homomorphism and we get the following important sufficient
condition for the vanishing of $\wt{MH}^{*,*}_{0}({\mcal X})$.
\begin{lemma}
\llabel{van0}
Let ${\cal X}$ be a pointed simplicial sheaf such that
$l\wt{H}^{*,*}({\cal X},\zz_{(l)})=0$. Then $\wt{MH}^{*,*}_{0}({\mcal
X})=0$.
\end{lemma}
\begin{proof}
Follows easily from the long exact sequences in the motivic cohomology
defined by the short exact sequence $0\sr \zz\stackrel{l}{\sr}\zz\sr
\zz/l\sr 0$.
\end{proof}
For a smooth variety $X$ over $k$ we denote by $s_d(T_X)$ the d-th
Milnor class of $X$ i.e. the characteristic class of the tangent
bundle $T_{X}$ which corresponds to the Newton polynomial in the Chern
classes (see \cite[Corollary 14.3]{Red} for a more careful
definition).
\begin{theorem}
\llabel{newth2} Let $Y$ be a smooth projective variety over $k$ such
that ther exists a map $X\sr Y$ where $X$ is a smooth projective
variety of dimension $l^m-1$ satisfying
$$deg(s_{l^m-1}(X))\ne 0(mod\,\,l^2)$$ 
Then $\wt{MH}^{p,q}_m(\wt{C}(Y),\zz/l)=0$ for all $p$ and $q$ (see
Appendix B. for the definition of $\wt{C}(Y)$).
\end{theorem}
\begin{proof}
For $m=0$ our condition means that $Y$ has a rational point over a
separable extension of $k$ of degree $n$ where $n\ne 0(mod\,\,
l^2)$. By Lemma \ref{van0a} this implies that
$l\wt{H}^{*,*}(\wt{C}(Y),\zz_{(l)})=0$ and by Lemma \ref{van0} we
conclude that $\wt{MH}^{*,*}_0(\wt{C}(Y))=0$.

Assume that $m>0$. Set $d=l^m-1$. Let $V$ be a vector bundle on $X$
and $f_V:T^{d+n}\sr Th_X(V)$ a map satisfying the conclusion of
Theorem \ref{newth8}. Consider the cofibration sequence
$$T^{d+n}\stackrel{f_V}{\sr} Th_X(V)\sr cone\stackrel{\psi}{\sr}
\Sigma^1_s T^{n+d}$$
corresponding to $f_V$. The long exact sequence of motivic cohomology
corresponding to this cofibration sequence shows that ther exists a
unique class $\alpha$ in $\tilde{H}^{2n,n}(cone)$ whose restriction to
$Th_X(V)$ is the Thom class $t_V$. Multiplication with this class
gives us a map
\begin{equation}
\llabel{prephi}
\wt{H}^{p,q}(\wt{C}(Y))\sr \wt{H}^{p+2n,q+n}(cone\wedge \wt{C}(Y))
\end{equation}
Consider the map $cone\wedge \wt{C}(Y)\sr (\Sigma^1_s
T^{n+d})\wedge\wt{C}(Y)$. We claim that this map is a weak
equivalence. Indeed, it is a part of a cofibration sequence and to
verify that it is a weak equivalence it is enough to check that
$Th_X(V)\wedge \wt{C}(Y)$ is contractible. This follows from the
cofibration sequence
$$(V-\{0\})_+\wedge \wt{C}(Y)\sr V_+\wedge \wt{C}(Y)\sr
Th_X(V)\wedge \wt{C}(Y)$$
and Lemma \ref{cc}. Since $cone\wedge \wt{C}(Y)\sr (\Sigma^1_s
T^{n+d})\wedge\wt{C}(Y)$ is a weak equivalence, (\ref{prephi}) defines
a map
$$\phi:\wt{H}^{p,q}(\wt{C}(Y))\sr \wt{H}^{p-2d-1,q-d}(\wt{C}(Y))$$
We claim that $\phi$ is something like a contracting homotopy for
the complex $(\wt{H}^{*,*}_n(\wt{C}(Y)), Q_m)$. More precisely we
have the following lemma which clearly imply the statement of the
theorem.
\begin{lemma}
\llabel{newl10}
There exists $c\in (\zz/l)^*$ such that for any $x\in
\wt{H}^{p,q}(\wt{C}(Y),\zz/l)$ one has
$$cx=\phi Q_m(x)-Q_m\phi(x)$$
\end{lemma}
\begin{proof}
Let $\gamma\in \wt{H}^{2n+2d+1,n+d}(cone)$ be the pull-back of the
tautological motivic cohomology class on $\Sigma^1_s T^{n+d}$ with
respect to $\psi$.  Since the map 
$$cone\wedge \wt{C}(Y)\sr
(\Sigma^1_s T^{n+d}) \wedge\wt{C}(Y)$$
is a weak equivalence it is sufficient to verify that there exists
$c\in (\zz/l)^*$ such that
$$c\gamma\wedge x=\alpha\wedge Q_m(x)-Q_m(\alpha\wedge x)$$
Let us show that for $i<m$ we have $Q_i(\alpha)=0$. Indeed, the
restriction of $Q_i(\alpha)$ to $Th_X(V)$ is zero by \cite[Theorem
14.2(1)]{Red} for any $i$. On the other hand for $i<m$ the restriction
map is a mono since there are no motivic cohomology of negative
weight. This fact together with \cite[Proposition 13.4]{Red} implies
that we have
$$Q_m(\alpha\wedge x)=\alpha\wedge Q_m(x) + Q_m(\alpha)\wedge x$$
and therefore
$$Q_m(\alpha)\wedge x=\alpha\wedge Q_m(x)-Q_m(\alpha\wedge x)$$
We now have two possibilities. If $Y$ has a point over an extension of
degree prime to $l$ then $\wt{H}^{p,q}(\wt{C}(X))=0$ and the statement of
our lemma obviously holds.  Assume that $Y$ has no points over
extensions of degree prime to $l$. Let us show that under this
assumption $Q_m(\alpha)=c\gamma$ for $c\in (\zz/l)^*$. Note first that
since $Y$ has no points over extensions of degree prime to $l$, $X$
does not have points over such extensions either. Since $Q_m(\alpha)$
restricts to zero on $Th_X(V)$ by \cite[Theorem 14.2(1)]{Red} it is
sufficient to show that $Q_m(\alpha)\ne 0$. By \cite[Proposition
13.6]{Red} we have $Q_m=\beta q_m\pm q_m\beta$. Since $\alpha$ is a
reduction of an integral class we have $Q_m(\alpha)=\beta q_m(\alpha)$
and to show that it is non-zero we have to check that $q_m(\alpha)$
can not be lifted to cohomology with $\zz/l^2$-coefficients. If it
could there would be a lifting $y$ of $q_m(t_V)$ to the cohomology
with $\zz/l^2$-coefficients such that $f_V^*(y)=0$. Our condition that
$X$ has no points over extensions of degree prime to $l$ implies, by
Lemma \ref{newl12}, that the value of $f_V^*(y)$ does not depend on
the choice of $y$. On the other hand by \cite[Corollary
14.3]{Red} we know that $q_m(t_V)$ is the reduction of the integral
class $s_{d}(X)$ which we assumed is non-zero mod $l^2$.
\end{proof}
\end{proof}
The following result provides us with a class of varieties $X$
satisfying the condition of Theorem \ref{newth2}.
\begin{proposition}
\llabel{newpr5}
Let $X$ be a smooth hypersurface of degree $d$ in ${\bf P}^n$. Then 
$$deg(s_{n-1}(X)):=deg(s_{n-1}(T_X))=d(n+1-d^{n-1})$$
In particular, if $X$ is a smooth hypersurface of degree $l$ in ${\bf
P}^{l^n}$ then 
$$deg(s_{l^n-1}(X))=l(mod\,\, l^2)$$
\end{proposition}
\begin{proof}
For the hypersurface given by a generic section of a vector bundle $L$
the normal bundle is canonically isomorphic to $L$. In particular we
have a short exact sequence of the form
$$0\sr T_X\sr i^*(T_{{\bf P}^n})\sr i^*({\cal O}(d))\sr 0$$
The tangent bundle on ${\bf P}^n$ fits into an exact sequence of the
form 
$$0\sr {\cal O}\sr {\cal O}(1)^{n+1}\sr T_{{\bf P}^n}\sr 0$$
Therefore, in $K_0(X)$ we have
$$T_X=i^*(T_{{\bf P}^n})-i^*({\cal O}(d))=(n+1)i^*({\cal O}(1))-{\cal
O}-i^*({\cal O}(d))$$
Since $s_{n-1}$ is an additive characteristic class which on line bundles
is given by $s_{n-1}(L)=e(L)^{n-1}$ we get for $n>1$
$$s_{n-1}(T_X)=(n+1-d^{n-1})i^*((-\sigma)^{n-1})$$
where $\sigma=e({\cal O}(-1))$. By Lemma \ref{newl7} we conclude that 
$$deg(s_{n-1}(X)):=deg(s_{n-1}(T_X))=d(n+1-d^{n-1})$$
\end{proof} 
Combining Theorem \ref{newth2} and Proposition \ref{newpr5} we get the
following result which is the only result of this section used for the
proof of Theorem \ref{maintheorem}.
\begin{cor}
\llabel{exitcor}
Let $Q$ be a smooth quadric in ${\bf P}^{2^n}$. Then
$$\wt{MH}^{*,*}_i(\wt{C}(Q),\zz/2)=0$$
for all $i\le n$.
\end{cor}

\subsection{Norm quadrics and their motives}
\label{third}\llabel{sec3}
The goal of this section is to prove Theorems \ref{Pmain} and
\ref{rostvan}. Everywhere in this section (except for Lemmas
\ref{sbs3}, \ref{sbs31} and \ref{cor3}) $k$ is a field of
characteristic $\ne 2$.

For elements
$a_1,\dots,a_n$ in $k^*$ let ${\langle}a_1,\dots,a_n{\rangle}$ be the quadratic form
$\sum a_i x_i^2$. One defines the Pfister form ${\langle\langle}a_1,\dots,a_n{\rangle\rangle}$ as
the tensor product
\begin{equation}
\llabel{nq}
{\langle\langle}a_1,\dots,a_n{\rangle\rangle}:={\langle}1,-a_1{\rangle}\oo\dots\oo{\langle}1,-a_n{\rangle}
\end{equation} 
Denote by $Q_{\uu{a}}=Q_{a_1,\dots,a_n}$ the projective quadric of
dimension $2^{n-1}-1$ given by the equation
${\langle\langle}a_1,\dots,a_{n-1}{\rangle\rangle}=a_nt^2$. For $n=2$ the rational points of the
affine part of this quadric ($t\ne 0$) can be identified with non-zero
elements $x$ of $E=k(\sqrt{a_1})$ such that $N_{E/k}(x)=a_2$. Because
of this interpretation the quadric given by (\ref{nq}) is called the
norm quadric associated with the sequence $(a_1,\dots,a_n)$. 

The following result is well known but we decided to include the proof
since it is crucial for our main theorem.
\begin{proposition}
\llabel{isspl0} The symbol $\uu{a}$ is divisible by $2$ in
$K_n^M(k(Q_{\uu{a}}))$.
\end{proposition}
\begin{proof}
We are going to show that if $Q_{\uu{a}}$ has a rational point over
$k$ then $\uu{a}$ is divisible by $2$ in $K_n^M(k)$. Since any variety
has a pojnt over its function field this implies that $\uu{a}$ is
divisible by 2 in the generic point of $Q_{\uu{a}}$. Let $P_{\uu{a}}$
denote the quadric given by the equation ${\langle\langle}a_1,\dots,a_n{\rangle\rangle}=0$.
\begin{lemma}
\llabel{l7}
For any $\uu{a}=(a_1,\dots,a_n)$ the following two conditions are
equivalent
\begin{enumerate}
\item $Q_{\uu{a}}$ has a $k$-rational point	
\item $P_{\uu{a}}$ has a $k$-rational point	
\end{enumerate}
\end{lemma}
\begin{proof}
The first condition implies the second one because the form
$${\langle\langle}a_1,\dots,a_{n-1}{\rangle\rangle}\oplus{\langle}-a_n{\rangle}$$
is a subform in ${\langle\langle}a_1,\dots,a_n{\rangle\rangle}$ and
therefore $Q_{\uu{a}}$ is a subvariety in $P_{\uu{a}}$. Assume that
the second condition holds. By \cite[Corollary 1.6]{Lam} it implies
that the form ${\langle\langle}a_1,\dots,a_n{\rangle\rangle}$ is
hyperbolic. Hence, for any rational point $p$ of $P_{\uu{a}}$ there
exists a linear subspace $H$ of dimension
$dim(P_{\uu{a}})/2=2^{n-1}-1$ which lies on $P_{\uu{a}}$ and passes
through $p$. The quadric $Q_{\uu{a}}$ is a section of $P_{\uu{a}}$ by
a linear subspace $L$ of codimension $2^{n-1}-1$ in ${\bf
P}^{2^n-1}$. The intersection of $H$ and $L$ is a rational point on
$Q_{\uu{a}}$.
\end{proof}
To prove Proposition \ref{isspl0} we proceed by
induction on $n$. Consider first the case $n=2$. Then $Q_{\uu{a}}$ is
given by the equation $x^2-a_1y^2=a_2z^2$. We may assume that it has a
point of the form $(x_0,y_0,1)$ (otherwise $a_1$ is a square root in $k$
and the statement is obvious). Then $a_2$ is the norm of the element
$x_0+\sqrt{a_1}y_0$ from $k(\sqrt{a_1})$ and thus the symbol $(a_1,a_2)$
is divisible by 2.

Suppose that the proposition is proved for sequences $(a_1,\dots,a_i)$
of length smaller than $n$. The quadric $Q_{\uu{a}}$ is given by the
equation 
$${\langle\langle}a_1,\dots,a_{n-1}{\rangle\rangle}=a_nt^2.$$
The form ${\langle\langle}a_1,\dots,a_{n-1}{\rangle\rangle}$ is of the
form ${\langle}1{\rangle}\oplus \mu_{\uu{a}}$. By induction we may
assume that our point $q\in Q_{\uu{a}}(k)$ belongs to the affine part
$t\ne 0$. Consider the plane $L$ generated by points $(1,0,\dots,0)$
and $q$. The restriction of the quadratic form
${\langle\langle}a_1,\dots,a_{n-1}{\rangle\rangle}$ to $L$ is of the
form ${\langle\langle}b{\rangle\rangle}$ for some $b$ (the idea is
that $L$ is a ``subfield'' in the vector space where
${\langle\langle}a_1,\dots,a_{n-1}{\rangle\rangle}$ lives). Consider
the field extension $k(\sqrt{b})$. The form
${\langle\langle}b{\rangle\rangle}$ and therefore the form
${\langle\langle}a_1,\dots,a_{n-1}{\rangle\rangle}$ represents zero
over $k(\sqrt{b})$ and thus by the inductive assumption
$(a_1,\dots,a_{n-1})=0$ in $K_{n-1}^M(k(\sqrt{b}))/2$.  On the other
hand by the construction ${\langle\langle}b{\rangle\rangle}$
represents $a_n$ and therefore we have $a_n\in Im
N_{k(\sqrt{b})/k}\subset k^*$ which proves the proposition.
\end{proof}
Let ${\cal X}_{\uu{a}}$ denote the simplicial scheme
$\check{C}(Q_{\uu{a}})$ (see Appendix B). The following result is
formulated in terms of the triangulated category of mixed motives
$DM_{-}^{eff}=DM^{eff}_{-}(k)$ introduced in \cite{H3new}.
\begin{theorem}
\llabel{Pmain} Let $k$ be a perfect field. Then, there exists a
distinguished triangle in $DM_{-}^{eff}$ of the form
\begin{equation}
\llabel{maindt}
M({\cal X}_{\uu{a}})(2^{n-1}-1)[2^n-2]\sr M_{\uu{a}}\sr M({\cal
X}_{\uu{a}}){\sr}M({\cal
X}_{\uu{a}})(2^{n-1}-1)[2^n-1]
\end{equation} 
\end{theorem}
\begin{proof}
Denote by $\cal C$ the category of Chow motives over $k$ and let $\zz\{n\}$ be
the n-th Tate motive in this category.  The proof of the theorem is
based on the following important result.
\begin{theorem}[Markus Rost]
\llabel{Rost7} There exists a direct summand $M_{\uu{a}}$ of
$Q_{\uu{a}}$ in $\cal C$ together with two morphisms
$$\psi^*:\zz\{2^{n-1}-1\}\sr M_{\uu{a}}$$
$$\psi_*:M_{\uu{a}}\sr \zz$$
such that
\begin{enumerate} 
\item the composition $\psi_*:Q_{\uu{a}}\sr M_{\uu{a}}\sr \zz$
is the morphism defined by the projection $Q_{\uu{a}}\sr Spec(k)$ 
\item for any field $F$ over $k$ where $Q_{\uu{a}}$ has a point the
pull-back of the sequence 
$$\zz\{2^{n-1}-1\}\sr M_{\uu{a}}\sr \zz$$
to $F$ is split-exact.
\end{enumerate}
\end{theorem}
\begin{proof}
See \cite{Rost2}, \cite{Rost2a}.
\end{proof}
The Friedlander-Lawson moving lemma for families of cycles shows that
for any $k$ there is a functor from the category of Chow motives over
$k$ to $DM_{-}^{eff}(k)$ (it is shown in \cite{H3new} that this
functor is a full embedding if $k$ is perfect). Therefore, the Rost
motive $M_{\uu{a}}$ can be also considered in $DM_{-}^{eff}$ where it
is a direct summand of the motive $M(Q_{\uu{a}})$ of the norm
quadric. To show that it fits into a distinguished triangle of the
form (\ref{maindt}) we need the following lemmas.

\begin{lemma}
\llabel{sbs3} Let $k$ be a perfect field, $M$ an object of
$DM^{eff}_-(k)$ and $X$ a smooth variety over $k$ such that for any
generic point $\eta$ of $X$ the pull-back of $M$ to the residue field
$k_{\eta}$ is zero. Then one has 
$$Hom_{DM}(M(X),M)=0$$
$$M\oo M(X)=0.$$
\end{lemma}
\begin{proof}
Since $k$ is perfect, the group $Hom(M(X),M)$ is the hypercohomology
group of $X$ with coefficients in the complex of sheaves with
transfers with homotopy invariant cohomology sheaves $K$ which
represents $M$. Our condition on $X$ and $M$ implies easily that the
cohomology sheaves $h_i(K)$ of $K$ vanish on the generic points of
schemes etale over $X$. Since $h_i(K)$ are homotopy invariant sheaves
with transfers we conclude by \cite[Corollary 4.19, p.116]{H2new} that
they vanish on all schemes etale over $X$. Therefore, $Hom(M(X),M)=0$.

To prove that $M\oo M(X)=0$ it is sufficient to show that the class of
objects $N$ such that
$$Hom(N,M\oo M(X)[m])=0$$
for all $m$, contains $M\oo M(X)$. Since this class is a localizing
subcategory\footnote{A subcategory in a triangulated category is
called localizing if it is closed under triangles, direct sums and
direct summands.}, it is sufficient to show that for any smooth $Y$
over $k$ and any $m$ one has
\begin{equation}
\llabel{sbs2}
Hom(M(Y\times X), M\oo M(X)[m])=0
\end{equation}
If $X$ and $M$ satisfy the conditions of the lemma then for any $Y$
and any $N$, $Y\times X$ and $M\oo N$ satify these
conditions. Therefore, (\ref{sbs2}) follows from the first assertion
of the lemma.
\end{proof}
\begin{lemma}
\llabel{tr1} Let $k$ be a perfect field. Then the sequence of
\begin{equation}
\llabel{sbs}
M(Q_{\uu{a}})(2^{n-1}-1)[2^n-2]\stackrel{Id\oo\psi^*}{\sr} M(Q_{\uu{a}})\oo M_{\uu{a}}\stackrel{Id\oo\psi_*}{\sr}
M(Q_{\uu{a}})
\end{equation}
is split-exact.
\end{lemma}
\begin{proof}
Observe first that $Id\oo \psi_*$ is a split epimorphism. Indeed,
Theorem \ref{Rost7}(1) implies that the morphism
$$M(Q_{\uu{a}})\sr M(Q_{\uu{a}})\oo M(Q_{\uu{a}})\sr M(Q_{\uu{a}})\oo
M_{\uu{a}}$$ 
where the first arrow is defined by the diagonal and the second by the
projection $M(Q_{\uu{a}})\sr M_{\uu{a}}$ is a section of $Id\oo
\psi_*$. It remains to show that (\ref{sbs}) extends to a
distinguished triangle. Let $cone$ be a cone of the morphism
$\psi^*:\zz(2^{n-1}-1)[2^n-2]\sr M_{\uu{a}}$. Since there are no
motivic cohomology of negative weight, the morphism
$\psi_*:M_{\uu{a}}\sr \zz$ factors through a morphism
$\phi:cone\sr\zz$. Let $cone'$ be the a cone of $\phi$. Standard
properties of triangles in triangulated categories imply that the
sequence (\ref{sbs}) extends to a distinguished triangle if and only
if $cone'\oo M(Q_{\uu{a}})=0$. It follows from Theorem \ref{Rost7}(b)
and Lemma \ref{sbs3}.
\end{proof}
\begin{lemma}
\llabel{sbs31} Let $X$ be a smooth scheme over $k$ and $M$ an
object of the localizing subcategory generated by $M(X)$. Then one has
\begin{enumerate}
\item the morphism $M\oo M(\check{C}(X))\sr M$ is an isomorphism
\item the homomorphism $Hom(M,\check{C}(X))\sr Hom(M,\zz)$ is an
isomorphism.
\end{enumerate}
\end{lemma}
\begin{proof}
It is clearly sufficient to prove the lemma for $M=M(X)$. In this case
the first statement follows immediately from Lemma \ref{cc} and the
fact that $M$ takes simplicial weak equivalences to isomorphisms. Let
$M(\wt{C}(X))$ be the cone of the morphism $M(\check{C}(X))\sr
\zz$. To prove the second statement we have to show that any morphism
$f:M(X)\sr M(\wt{C}(X))$ is zero. The morphism $f$ can be written as
the composition
$$M(X)\stackrel{M(\Delta)}{\sr}M(X)\oo M(X)\stackrel{f\oo
Id}{\sr}M(\wt{C}(X))\oo M(X)\sr  M(\wt{C}(X))$$
which is zero because the first part of the lemma implies that
$$M(\wt{C}(X))\oo M(X)=0.$$
\end{proof}
{\bf (Proof of Theorem \ref{Pmain} continues)} The morphism
$\psi_*:M_{\uu{a}}\sr \zz$ has a canonical lifting to a morphism
$\tilde{\psi}_*:M_{\uu{a}}\sr M({\cal X}_{\uu{a}})$ by Lemma
\ref{sbs31}(2). Together with the composition
$$\tilde{\psi}^*:M({\cal X}_{\uu{a}})(2^{n-1}-1)[2^n-2]\sr
\zz(2^{n-1}-1)[2^n-2]\stackrel{\psi^*}{\sr}M_{\uu{a}}$$
this lifting gives us a sequence of morphisms
$$M({\cal X}_{\uu{a}})(2^{n-1}-1)[2^n-2]\stackrel{\tilde{\psi}^*}{\sr}
M_{\uu{a}}\stackrel{\tilde{\psi}_*}{\sr} M({\cal 
X}_{\uu{a}})$$
The composition $\tilde{\psi}_*\circ\tilde{\psi}^*$ is zero by Lemma
\ref{sbs31}(2) and the fact that 
$$Hom(\zz(2^{n-1}-1)[2^n-2],\zz)=0$$
Let $cone$ be a cone of $\tilde{\psi}^*$. The morphism
$\tilde{\psi}_*$ factors through a morphism $\phi:cone\sr M({\cal
X}_{\uu{a}})$ and we have to show that $\phi$ is an isomorphism. The
category $\cal C$ of objects $N$ such that $\phi\oo Id_N$ is an
isomorphism is a localizing subcategory. By Lemma \ref{sbs31}(1), the
morphism $M({\cal X}_{\uu{a}})\oo M(Q_{\uu{a}})\sr M(Q_{\uu{a}})$ is
an isomorphism and we conclude by Lemma \ref{tr1} that $\cal C$
contains $M(Q_{\uu{a}})$ and therefore it contains $M({\cal
X}_{\uu{a}})$. On the other hand we have a commutative diagram
$$
\begin{array}{ccc}
cone\oo M({\cal X}_{\uu{a}})&\sr& M({\cal
X}_{\uu{a}})\oo M({\cal X}_{\uu{a}})\\
\downarrow&&\downarrow\\
cone&\sr&M({\cal X}_{\uu{a}})
\end{array}
$$
with both vertical arrows are isomorphisms by Lemma
\ref{sbs31}(1). This finishes the proof of Theorem \ref{Pmain}.
\end{proof}
\begin{remark}\rm
It is essential to use the category $DM_{-}^{eff}$ in Theorem
\ref{Pmain} because the triangle (\ref{maindt}) can not be lifted to
the motivic homotopy category (stable or unstable). One can see this
for $n=2$ using the fact that the map ${\mcal X}_{\uu{a}}\sr Spec(k)$
defines an isomorphism on algebraic K-theory but
$K^{*,*}(Q_{\uu{a}})\ne K^{*,*}(k)\oplus K^{*,*}(k)$.
\end{remark}
\begin{theorem}
\llabel{rostvan} $H^{2^{n}-1,2^{{n}-1}}({\cal
X}_{\uu{a}},\zz_{(2)})=0$
\end{theorem}
\begin{proof}
Since $char(k)\ne 2$ any purely inseparable extension of $k$ is of odd
degree and the transfer argument shows that it is sufficient to
consider the case of a perfect $k$. Denote by $\uu{K}_n^M$ the sheaf
on $Sm/k$ such that for a connected smooth scheme $X$ over $k$ the
group $\uu{K}^M_n(X)$ is the subgroup in the n-th Milnor K-group of
the function field of $X$ which consists of elements $u$ such that all
residues of $u$ in points of codimension $1$ are zero. The proof is
based on the following result.
\begin{theorem}[Markus Rost]
\llabel{l9} The natural homomorphism
\begin{equation}
\llabel{mainmono}
H^{2^{n-1}-1}(Q_{\uu{a}},\uu{K}_{2^{n-1}}^M)\sr k^*
\end{equation}
is a monomorphism.
\end{theorem}
\begin{proof}
See \cite{Rost1}.
\end{proof}
The following lemma shows that the cohomology group on the left hand
side of (\ref{mainmono}) can be replaced by a motivic cohomology
group.
\begin{lemma}
\llabel{cor3} Let $X$ be a smooth scheme over a field $k$. Then
for any $p,q$ there is a canonical homomorphism $H^{p,q}(X,\zz)\sr
H^{p-q}(X,\uu{K}^M_q)$ which is an isomorphism if $p\ge q+dim(X)$.
\end{lemma}
\begin{proof}
Considering $X$ as a limit of smooth schemes over the subfield of
constants of $k$ we may assume that $k$ is perfect. Let
$h^i=\uu{H}^i(\zz(q))$ denote the cohomology sheaves of the complex
$\zz(q)$. Since $h^i=0$ for $i>q$ the standard spectral sequence which
goes from cohomology with coefficients in $h^i$ and converges to the
hypercohomology with coefficients in $\zz(q)$ shows that there is a
canonical homomorphism
\begin{equation}
\llabel{homom}
H^{p,q}(X,\zz)\sr H^{p-q}(X,h^q)
\end{equation}
The same spectral sequence implies that the kernel and cokernel of
this homomorphism are bulit out of groups $H^{p-i}(X,
h^i)$ and $H^{p-i+1}(X, h^i)$
respectively, where $i<q$.  Since $p\ge q+dim(X)$ we get
$p-i>p-q\ge dim(X)$ and the cohomological dimension theorem for the
Nisnevich topology implies that these groups are zero.

It remains to show that $h^q=\uu{K}_q^M$. Since
$h^q$ is a homotopy invariant sheaf with transfers for
any smooth connected $X$ the restriction homomorphism
$$H^0(X,h^q)\sr H^0(Spec(k(X)),h^q)$$
is a monomorphism (\cite[Corollary 4.18, p.116]{H2new}).
It was shown in \cite[Prop. 4.1]{SusVoe3} that for any field $k$ one
has canonical isomorphisms $H^{q,q}(k,\zz)=K_q^M(k)$. In particular
for any $X$ the group $H^0(X,h^q)$ is a subgroup in
$K_q^M(k(X))$ and one verifies easily that it coincides with the
subgroup $\uu{K}_q^M(X)$ of elements with zero residues.
\end{proof}
By Theorem \ref{Pmain} and the suspension isomorphism (see
\cite[Theorem 2.4]{Red}) we have an exact sequence
\begin{equation}
\llabel{longexold}
H^{0,1}({\cal X}_{\uu{a}},\zz)\sr H^{2^n-1,2^{n-1}}({\cal
X}_{\uu{a}},\zz)\sr H^{2^n-1,2^{n-1}}(M_{\uu{a}},\zz)\sr H^{1,1}({\cal
X}_{\uu{a}},\zz)
\end{equation}
Since the dimension of $Q_{\uu{a}}$ equals $2^{n-1}-1$, the left hand
side group in (\ref{mainmono}) is isomorphic to the group $H^{2^{n}-1,
2^{n-1}}(Q_{\uu{a}},\zz)$ by Lemma \ref{cor3} and we obtain a
natural monomorphism
$$H^{2^{n}-1, 2^{n-1}}(Q_{\uu{a}},\zz)\sr k^*.$$
Let $\bar{M}_{\uu{a}}$ be the lifting of $M_{\uu{a}}$ to the algebraic
closure of $k$.  Since $M_{\uu{a}}$ is a direct summand of
$M(Q_{\uu{a}})$ we conclude that the map 
$$H^{2^n-1,2^{n-1}}(M_{\uu{a}},\zz)\sr
H^{2^n-1,2^{n-1}}(\bar{M}_{\uu{a}},\zz)$$
is injective. Since
$$H^{2^n-1,2^{n-1}}({\cal
X}_{\uu{a}}\times
Spec(\bar{k}),\zz)=H^{2^n-1,2^{n-1}}(Spec(\bar{k}),\zz)=0$$
we conclude that the second arrow in (\ref{longexold}) is zero. On the
other hand since $\zz(1)={\bf G}_m[-1]$ we have
$$H^{0,1}({\cal X}_{\uu{a}},\zz)={H}^{-1}({\cal X}_{\uu{a}},{\bf
G}_m)=0$$
\end{proof}

\subsection{Computations with Galois cohomology}\llabel{sec4}
In this section we are only concerned with classical objects - 
Milnor K-theory and etale cohomology. More general motivic cohomology
does not appear here. The only result of this section which we will
directly use below is Theorem \ref{l4}. One may observe that in
the case of $\zz/2$-coefficients it can be proved in a much easier
way, but we decided to include the case of general $l$ for possible
future use.
\begin{definition}
\llabel{BKd} We say that $BK({w},l)$ holds if for any field $k$ of
characteristic $\ne l$ and any $q\le {w}$ one has:
\begin{enumerate}
\item the norm residue homomorphism $K_q^M(k)/l\sr H^q_{et}(k,{\bf
\mu}_l^{\oo q})$ is an isomorphism
\item for any cyclic extension $E/k$ of degree $l$ the sequence
$$K_{q}^M(E)\stackrel{1-\sigma}{\sr}K_{q}^M(E)\stackrel{N_{E/k}}{\sr}
K_q^M(k)$$
where $\sigma$ is a generator of $Gal(E/k)$ is exact.
\end{enumerate}
\end{definition}
\begin{proposition}
\llabel{l0} Let $k$ be a field of characteristic $\ne l$ which has no
extensions of degree prime to $l$. Assume that $BK({w},l)$ holds. Then
for any cyclic extension $E/k$ of $k$ of degree $l$ there is an exact
sequence of the form
$$H^{w}_{et}(E,\zz/l)\stackrel{N_{E/k}}{\sr}H^{w}_{et}(k,\zz/l)
\stackrel{-\wedge[a]}{\sr}H^{{w}+1}_{et}(k,\zz/l)\sr H^{{w}+1}_{et}(E,\zz/l)$$
where $[a]\in H^1_{et}(k,\zz/l)$ is the class which corresponds to
$E/k$.  
\end{proposition}
\begin{proof}
In order to prove the proposition we will need to do some preliminary
computations.  Fix an algebraic closure $\bar{k}$ of $k$. Since $k$
has no extensions of degree prime to $l$ there exists a primitive root
of unit $\xi\in k$ of degree $l$. Let $E\subset \bar{k}$ be a cyclic
extension of $k$ of degree $l$. We have $E=k(b)$ where $b^l=a$ for an
element $a$ in $k^*$.  Denote by $\sigma_{\xi}$ the generator of the
Galois group $G_b=Gal(E/k)$ which acts on $b$ by multiplication by
$\xi$ and by $[a]_{\xi}$ the class in $H^1_{et}(k,\zz/l)$ which
corresponds to the homomorphism $Gal(\bar{k}/k)\sr G\sr \zz/l$ which
takes $\sigma_{\xi}$ to the canonical generator of $\zz/l$ (one can
easily see that this class is determined by $a$ and $\xi$ and does
notdepend on $b$).

Let $p:Spec(E)\sr Spec(k)$ be the projection. Consider the etale sheaf
$F=p_*(\zz/l)$. The group $G$ acts on $F$ in the natural way. Denote by
$F_i$ the kernel of the 
homomorphism $(1-\sigma)^i:F\sr F$. One can verify easily that
$F_i=Im(1-\sigma)^{l-i}$ and that as a $\zz/l[Gal(\bar{k}/k)]$-module
$F_i$ has dimension $i$. In particular we have $F=F_l$. Note that the
extension 
$$0\sr \zz/l\sr F_2\sr \zz/l\sr 0$$
represents the element $[a]_{\xi}$ in
$H^1_{et}(k,\zz/l)=Ext^1_{\zz/l[Gal(\bar{k}/k)]}(\zz/l,\zz/l)$. Let
$\alpha_i$ be the element in
$H^2(k,\zz/l)=Ext^2_{\zz/l[Gal(\bar{k}/k)]}(\zz/l,\zz/l)$ defined by
the exact sequence
$$0\sr\zz/l\sr F_i\stackrel{u_i}{\sr}F_i\sr\zz/l\sr 0$$
where $u_i=1-\sigma$ and $Im(u_i)=F_{i-1}$.
\begin{lemma}
\llabel{comp1}
One has $\alpha_l = c([a]_{\xi}\wedge[\xi]_{\xi})$ where
$c$ is an invertible element of $\zz/l$ and $\alpha_i=0$ for $i<l$. 
\end{lemma}
\begin{proof}
The fact that $\alpha_i=0$ for $i<l$ follows from the commutativity of
the diagram
$$
\begin{array}{ccccccccc}
0\sr&\zz/l&\sr&F_{i+1}&\sr&F_{i+1}&\sr&\zz/l&\sr 0\\
&0\downarrow&&\downarrow&&\downarrow&&\downarrow Id&\\
0\sr&\zz/l&\sr&F_{i}&\sr&F_{i}&\sr&\zz/l&\sr 0
\end{array}
$$

To compute $\alpha_l$ note first that since the action of
$Gal(\bar{k}/k)$ on $F=F_l$ factors through $G=Gal(E/k)$ it comes from a
well defined element in $H^2(G,\zz/l)$. This element is not
zero for trivial reasons. On the other hand the group 
$$H^2(G,\zz/l)=H^2(\zz/l,\zz/l)=\zz/l$$
is generated by the element $\beta(\gamma)$ where $\gamma$ is the
canonical generator of $H^1(G,\zz/l)$ and $\beta$ is the Bockstein
homomorphism. Thus we conclude that up to multiplication by an
invertible element of $\zz/l$ our class $\alpha_l$ equals
$\beta([a]_{\xi})$. It remains to show that
$\beta([a]_{\xi})=c[a]_{\xi}\wedge[\xi]_{\xi}$ which 
 follows by simple explicite computations from the fact that
$[a]_{\xi}$ has a lifting to an element of
$H^1_{et}(Spec(k),\mu_{l^2})$. 
\end{proof}
\begin{lemma}
\llabel{comp2} Assume that $BK({w},l)$ holds. Then for all fields $k$ of
characteristic $\ne l$, all $q\le {w}$ and all $i=1,\dots,l-1$ one has:
\begin{enumerate}
\item The sequence $H^q_{et}(k,\zz/l)\oplus H^q_{et}(k,F_{i+1})\sr
H^m_{et}(k,F_{i+1})\sr H^q_{et}(k,\zz/l)$ 
where the first homomorphism is given on the second summand by
$1-\sigma$ is exact.
\item The homomorphisms $\nu_{q,i}:H^q_{et}(k,\zz/l)\oplus
H^m_{et}(k,F_{i+1})\sr H^q_{et}(k, F_i)$ 
given by the canonical morphisms $\zz/l\sr F_{i}$, $F_{i+1}\sr
F_i$  are surjective.
\end{enumerate}
\end{lemma}
\begin{proof}
We proceed by induction on $i$. Consider first the case $i=l-1$. The
first statement follows immediately the assumption that $BK({w},l)$ holds. 

Let us prove the second one. The
image of $H^q_{et}(k, F_{l-1})$ in 
$$H^q_{et}(k,F_l)=H^q_{et}(E,\zz/l)$$
coincides with the kernel of the norm homomorphism
$$H^q_{et}(E,\zz/l)\sr H^q_{et}(k,\zz/l)$$ 
The first statement implies then that $H^q_{et}(k,\zz/l)\oplus
H^q_{et}(k,F_{l})$ maps surjectively to this image. It is sufficient
therefore to show that an element $\gamma\in H^q_{et}(k, F_{l-1})$ which
goes to zero in $H^q_{et}(k,F_l)$ belongs to the image of
$\nu_{q,l-1}$. Any such element is a composition of a cohomology class
in $H^{q-1}_{et}(k,\zz/l)$ with the canonical  extension
$$0\sr F_{l-1}\sr F_l\sr \zz/l\sr 0$$
Thus we may assume that $q=1$ and $\gamma$ is the element which
corresponds to this extension. Let  $\delta$ be the image of
$c[\xi]_{\xi}$ (where $c$ is as in Lemma \ref{comp1}) under the
homomorphism $H^1_{et}(k,\zz/l)\sr H^1_{et}(k,F_{l-1})$. The composition
$$H^1_{et}(k,\zz/l)\sr H^1_{et}(k,F_{l-1})\sr H^2_{et}(k,\zz/l)$$
where the later homomorphism corresponds to the extension
$$0\sr \zz/l\sr F_l\sr F_{l-1}\sr 0$$
equals to multiplication by $[a]_{\xi}$. We conclude now by Lemma
\ref{comp1} that the image of $\gamma-\delta$ in $H^2(k,\zz/l)$ is
zero. Then it lifts to $H^1_{et}(k,F_l)$ which proves our Lemma in the
case $i=l-1$. 

Suppose that the lemma is proved for all $i>j$. Let us show that it
holds for $i=j$. The first statement follows immediately from the
inductive assumption and the commutativity of the diagram
$$
\begin{array}{ccc}
F_{j+2}&\stackrel{1-\sigma}{\sr}&F_{j+2}\\
\downarrow&&\downarrow\\
F_{j+1}&\stackrel{1-\sigma}{\sr}&F_{j+1}
\end{array}
$$

The proof of the second one is now similar to the case $i=l-1$ with a
simplification due to the fact that $\alpha_i=0$ for $i<l$ (Lemma
\ref{comp1}).
\end{proof}
The statement of the proposition follows immediately from Lemma \ref{comp2}.
\end{proof}
\begin{remark}
\rm 
For $l=2$ Proposition \ref{l0} is a trivial corollary of the exactness
of the sequence $0\sr \zz/2\sr F_2\sr \zz/2\sr 0$. In particular it
holds without the $BK({w},l)$ assumption and not only in the
context of Galois cohomology but for cohomology of any
(pro-)finite group. 
For $l>2$ this is not true anymore which one can see
considering cohomology of $\zz/l$.
\end{remark}
\begin{lemma}
\llabel{l.5} Assume that $BK({w},l)$ holds and let $k$ be a field of
characteristic $\ne l$ which has no extensions of degree prime to
$l$. Let further $E/k$ be a cyclic extension of degree $l$ such that
the norm homomorphism $K_{w}^M(E)\sr K^M_{w}(k)$ is surjective. Then the
sequence
$$K_{{w}+1}^M(E)\stackrel{1-\sigma}{\sr}K_{{w}+1}^M(E)
\stackrel{N_{E/k}}{\sr}K_{{w}+1}^M(k)$$ 
where $\sigma$ is a generator for $Gal(E/k)$ is exact.
\end{lemma}
\begin{proof}
It is essentially a version of the proof given in \cite{Suslin4} for
${w}=2$ and
in \cite{MS2} for ${w}=3$. Let
us define a homomorphism 
$$\phi:K_{{w}+1}^M(k)\sr
K_{{w}+1}^M(E)/(Im(1-\sigma))$$
 as 
follows. Let $a$ be an  element in $K_{{w}+1}^M(k)$ of the form
$(a_0,\dots,a_{w})$ and let $b$ be an element in $K_{w}^M(E)$ such that 
$$N_{E/k}(b)=(a_0,\dots,a_{{w}-1})$$ We set $\phi(a)=b\wedge
a_n$. Since $BK({w},l)$ holds the element $\phi(a)$ does not depend on
the choice of $b$ and one can easily see that $\phi$ is a homomorphism
from $(k^*)^{\oo {w}+1}$ to $K_{{w}+1}^M(E)/(Im(1-\sigma))$. To show
that it is a homomorphism from $K_{w}^M(k)$ it is sufficient to verify
that $\phi$ takes an element of the form $(a_0,\dots,a_{w})$ such that
say $a_0+a_{w}=1$ to zero. Let $b$ be a preimage of
$(a_0,\dots,a_{{w}-1})$ in $K_{w}^M(k)$. We have to show that
$(b,a_{w})\in (1-\sigma)K_{{w}+1}^M(E)$. Assume first that $a_0$ is
not in $(k^*)^l$ and let $c$ be an element in $\bar{k}^*$ such that
$c^l=a_0$. Let further $F=k(c)$. Then by $BK({w},l)$ one has
$$b\wedge a_{w}=b\wedge(1-a_0)=N_{EF/E}(b_{EL}\wedge(1-c))=$$
$$=N_{EF/E}((b-(c,a_1,\dots,a_{w}-1))\wedge (1-c))\in
(1-\sigma)K_{{w}+1}^M(E)$$
since $N_{EF/F}(b-(c,a_1,\dots,a_{w}-1))=0$. The proof for the case when
$a_0\in (k^*)^l$ is similar. 

Clearly $\phi$ is a section for the obvious morphism
$$K_{{w}+1}^M(E)/(Im(1-\sigma))\sr K_{{w}+1}^M(k)$$
It remains to show that it is surjective. It follows immediately from
the fact that under our assumption on $k$ the group $K_{{w}+1}^M(E)$
is generated by symbols of the form $(b,a_1,\dots,a_{w})$ where $b\in
E^*$ and $a_1,\dots,a_{w}\in k^*$ (see \cite{BT}).
\end{proof}
\begin{lemma}
\llabel{l1} Assume that $BK({w},l)$ holds and let $k$ be a field of
characteristic $\ne l$ which has no extensions of degree prime to
$l$. Then the following two conditions are equivalent
\begin{enumerate}
\item $K_{{w}+1}^M(k)=l K_{{w}+1}^M(k)$
\item for any cyclic extension $E/k$ the norm homomorphism
$$K_{w}^M(E)\sr K_{w}^M(k)$$
is surjective.
\end{enumerate}
\end{lemma}
\begin{proof}
The $2\Rightarrow 1$ part follows from the projection formula. Since
$BK({w},l)$ holds we conclude that there is a commutative square with
surjective horizontal arrows of the form
$$
\begin{array}{ccc}
K_{w}^M(E)&\sr&H^{w}(E,\zz/l)\\
\downarrow&&\downarrow\\
K_{w}^M(k)&\sr&H^{w}(k,\zz/l)
\end{array}
$$
and thus the cokernel of the left vertical arrow is the same as the
cokernel of the right one. By Proposition \ref{l0} it gives us an exact
sequence
$$K_{w}^M(E)\sr K_{w}^M(k)\sr H^{w+1}(k,\zz/l)$$
and since the last arrow clearly factors through $K_{w+1}^M(k)/l$ it is
zero. Lemma is proved.
\end{proof}
\begin{lemma}
\label{l3}
Assume that $BK({w},l)$ holds and let $k$ be a field of characteristic
not equal to $l$ which has no extensions of degree prime to
$l$. Assume further that $K_{{w}+1}^M(k)=l K_{{w}+1}^M(k)$.  Then for
any finite extension $E/k$ one has $K_{{w}+1}^M(E)=l K_{{w}+1}^M(E)$.
\end{lemma}
\begin{proof}
This proof is a variant of the proof given in \cite{Suslin4} for
${w}=1$. Since $k$ has no extensions of degree prime to $l$ it is
separable and its Galois group is an l-group. Therefore it is
sufficient to prove the lemma in the case of a cyclic extension $E/k$
of degree $l$. By Proposition \ref{l0} and Lemma \ref{l.5} we have an
exact sequence
$$K_{{w}+1}^M(E)\stackrel{1-\sigma}{\sr}K_{{w}+1}^M(E)\stackrel{N_{E/k}}{\sr}K_{{w}+1}^M(k).$$
Let $\alpha$ be an element in $K_{{w}+1}^M(E)$ and let $\beta\in
K_{{w}+1}^M(k)$ be an element such that $N_{E/k}(\alpha)=l\beta$. Then
$N_{E/k}(\alpha-\beta_E)=0$ and we conclude that the endomorphism
$$1-\sigma:K_{{w}+1}^M(E)/l\sr K_{{w}+1}^M(E)/l$$
is surjective. Since $(1-\sigma)^l=0$ this implies that $K_{{w}+1}^M(E)/l=0$.
\end{proof}
\begin{theorem}
\llabel{l4} Assume that $BK({w},l)$ holds and let $k$ be a field of
characteristic not equal to $l$ which has no extensions of degree
prime to $l$ such that $K_{{w}+1}^M(k)=l K_{{w}+1}^M(k)$.  Then
$H^{{w}+1}_{et}(k,\zz/l)=0$.
\end{theorem}
\begin{proof}
Let $\alpha$ be an element of $H^{{w}+1}_{et}(k,\zz/l)$. We
have to show that $\alpha=0$. By Lemma \ref{l3} and obvious induction we
may assume that 
$\alpha$ vanishes on a cyclic extension of $k$. Then by Proposition
\ref{l0} it is of the form $\alpha_0\wedge a$ where $a\in
H^1_{et}(k,\zz/l)$ is the element which represents our cyclic
extension. Thus since $BK({w},l)$ holds it belongs to the image of the
homomorphism $K_{{w}+1}^M(k)/l\sr H^{{w}+1}_{et}(k,\zz/l)$ and therefore is zero.
\end{proof}
\subsection{Beilinson-Lichtenbaum conjectures}
\label{first}\llabel{sec5}
For a smooth variety $X$ and an abelian group $A$ define the {\em
Lichtenbaum motivic cohomology groups} of $X$ with
coefficients in $A$ as the hypercohomology groups
$$H^{p,q}_L(X,A):={\bf H}^p_{et}(X,A\oo\zz(q))$$ 
The following fundamental conjecture is due to Alexander Beilinson and
Stephen Lichtenbaum (see \cite{pairing},\cite{Licht2}).
\begin{conjecture}
\llabel{hilb90}
Let $k$ be a field. Then for any $n\ge 0$ one has 
$$H^{n+1,n}_L(Spec(k),\zz)=0$$
\end{conjecture}
For $n=0$ we have $\zz(0)=\zz$ and this conjecture follows from the
simple fact that $H^1(k,\zz)=0$. For $n=1$ we have $\zz(1)={\bf
G}_m[-1]$ and the conjecture is equivalent to the cohomological form
of the Hilbert 90 Theorem. Because of this fact Conjecture
\ref{hilb90} is called the generalized Hilbert 90 Conjecture. 

The standard proofs of Conjecture \ref{hilb90} for $n=0$ and $n=1$
work integrally. In the proofs of all the known cases of Conjecture
\ref{hilb90} for $n>1$ one considers the vanishing of the groups
$H^{n+1,n}_L(Spec(k),\zz_{(l)})=0$ for different primes $l$
separately.
\begin{definition}
\llabel{h90} We say that $H90(n,l)$ holds if for any $k$ one has
$$H^{n+1,n}_L(Spec(k),\zz_{(l)})=0.$$
\end{definition}
The following result is proved in \cite[Theorem
8.6]{GL2}.\footnote{The authors of \cite{GL2} define motivic
cohomology through the higher Chow groups. Their definition is
equivalent to the one used here by the comparison theorem of
\cite{Suslin3new}, \cite{FS} and \cite{comparison}.}
\begin{theorem}[Geisser-Levine]
\llabel{gl} Let $S$ be the local ring of a smooth scheme over a field
of characteristic $p>0$. Then for any $n\ge 0$ one has
$$H^{n+1,n}_L(Spec(k),\zz_{(p)})=0.$$
\end{theorem}
Let $\pi:(Sm/k)_{et}\sr (Sm/k)_{Nis}$ be the obvious morphism of
sites. Consider the complex ${\bf
R}\pi_*(\pi^*(\zz(q)))$ of Nisnevich sheaves with transfers on $Sm/k$.
We have
$$H^{p,q}_L(X,\zz)={\bf H}^{p}_{Nis}(X,{\bf R}\pi_*(\pi^*(\zz(q)))).$$
Let $L(q)$ be the canonical truncation of the complex ${\bf
R}\pi_*(\pi^*(\zz(q)))$ at level $q+1$ i.e. $L(q)$ is the subcomplex
of sheaves in ${\bf R}\pi_*(\pi^*(\zz(q)))$ whose cohomology sheaves
$\uu{H}^i(L(q))$ are the same as for ${\bf R}\pi_*(\pi^*(\zz(q)))$ for
$i\le q+1$ and zero for $i>q+1$. Since $\uu{H}^i(\zz(q))=0$ for $i>q$
the canonical morphism $\zz(q)\sr {\bf R}\pi_*(\pi^*(\zz(q)))$ factors
through $L(q)$.  Let $K(q)$ be the complex of sheaves with transfers
on $(Sm/k)_{Nis}$ defined by the distinguished triangle $\zz(q)\sr
L(q)\sr K(q) \sr \zz(q)[1]$.
\begin{theorem}
\llabel{sv} Assume that $H90({w},l)$ holds. Then for any $k$ such that
$char(k)\ne l$ the complex $K(w)\oo\zz_{(l)}$ is quasi-isomorphic to
zero.
\end{theorem}
\begin{proof}
We will use the following two lemmas.
\begin{lemma}
\llabel{hi} Let $k$ be a field of characteristic not equal to
$l$. Then for any $q\ge 0$ and any $\zz_{(l)}$-module $A$, the complex
${\bf R}\pi_*(\pi^*(A\oo\zz(q)))$ has homotopy invariant cohomology sheaves.
\end{lemma}
\begin{proof}
Since the sheaf associated with a homotopy invariant presheaf with
transfers is homotopy invariant (\cite[Propositions 4.26,
p. 118]{H2new} and \cite[Propositions 5.5, p. 128]{H2new}) it is
sufficient to show that the functors $H^{p,q}_L(-,A)$ are homotopy
invariant i.e. that for any smooth $U$ the homomorphism
$$H^{p,q}_L(U\times\af,A)\sr H^{p,q}(U,A)$$ given by the restriction
to $U\times\{0\}$ is an isomorphism. Consider the universal
coefficients long exact sequence relating Lichtenbaum motivic
cohomology with coefficients in $A$ to Lichtenbaum motivic cohomology
with coefficients in $A\oo\qq$ and $A\oo\qq/\zz_{(l)}$. The cohomology
with $A\oo\qq$-coefficients is homotopy invariant by Lemma 
\ref{ratLnew}. The cohomology with $A\oo\qq/\zz_{(l)}$-coefficients is
isomorphic to the etale cohomology by Theorem \ref{sv0} and therefore
homotopy invariant as well. Our claim follows now from the five lemma.
\end{proof}
\begin{lemma}
\llabel{ratLnew}
The canonical homomorphisms 
$$H^{p,q}(-,\qq)\sr H^{p,q}_L(-,\qq)$$
are isomorphisms. 
\end{lemma}
\begin{proof}
It is sufficient to show that for any $i,j\in\zz$ and any smooth
scheme $X$ one has
$$H^i_{et}(X,\uu{H}^j(\qq(q))_{et})=H^i_{Nis}(X,\uu{H}^j(\qq(q))).$$
Since the $\uu{H}^j(\zz(q))$ are sheaves with transfers it is a particular
case of \cite[Propositions 5.24, 5.27, p. 135]{H2new}.
\end{proof}
Lemma \ref{hi} implies that $\uu{H}^{n+1}(K({w})\oo\zz_{(l)})$ is a
homotopy invariant sheaf with transfers and by the assumption that
$H90({w},l)$ holds we know that it vanishes over fields. Therefore, by
\cite[Corollary 4.18, p. 116]{H2new} we conclude that
$\uu{H}^{n+1}(K({w})\oo\zz_{(l)})=0$.

In order to show that $K({w})\oo\zz_{(l)}$ is quasi-isomorphic to zero
it remains to verify that for any smooth scheme $X$ over $k$ and any
$p\le {w}$ the homomorphism $H^{p,{w}}(X,\zz_{(l)})\sr
H^{p,{w}}_L(X,\zz_{(l)})$ is an isomorphism. Lemma \ref{ratLnew} and
the universal coefficients long exact sequence imply that it is
sufficient to verify that the homomorphisms
\begin{equation}
\llabel{inques}
H^{p,{w}}(X,\qq/\zz_{(l)})\sr H^{p,{w}}_L(X,\qq/\zz_{(l)})
\end{equation} 
are isomorphisms for $p\le {w}$. Since we assume $H90({w},l)$, we
conclude that the group $H^{p,{w}}_L(X,\qq/\zz_{(l)})$ is infinitely
divisible and using the trick of \cite[Theorem 11.4]{SusVoe3} we
conclude that the map (\ref{inques}) is surjective for $p=w$ and
$X=Spec(F)$ where $F$ is a field. By the analog of \cite[Theorem
1.1]{GL} for $\qq/\zz_{(l)}$-coefficients and the comparison between
the higher Chow groups and the etale cohomology we conclude that
(\ref{inques}) is an isomorphism.
\end{proof}
\begin{cor}
\llabel{svnew} Assume that $H90({w},l)$ holds. Then for any field $k$
and any smooth simplicial scheme $\cal X$ over $k$ one has
\begin{enumerate}
\item the homomorphisms
$$H^{p,q}({\cal X},\zz_{(l)})\sr H^{p,q}_L({\cal X},\zz_{(l)})$$
are isomorphisms for $p-1\le q\le {w}$ and monomorphisms for $p=q+2$ and
$q\le {w}$ 
\item the homomorphisms
$$H^{p,q}({\cal X},\zz/l)\sr H^{p,q}_L({\cal X},\zz/l)$$
are isomorphisms for $p\le q\le {w}$ and monomorphisms for $p=q+1$ and
$q\le {w}$ 
\end{enumerate}
\end{cor}

For $n$ prime to the characteristic of the base field, Lichtenbaum
motivic cohomology groups with $\zz/n$ coefficients are closely related
to the ``usual'' etale cohomology.  
\begin{theorem}
\llabel{sv0}
Let $k$ be a field and $n$ be an integer prime to characteristic of
$k$. Denote by ${\bf \mu}_n$ the etale sheaf of n-th roots of unit on
$Sm/k$ and let ${\bf \mu}_n^{\oo q}$ be the n-th tensor power of ${\bf
\mu}_n$ in the category of $\zz/n$-modules. Then there is a canonical
isomorphism $H^{p,q}_{L}(-,\zz/n)=H^{p}_{et}(-,{\bf \mu}_n^{\oo q})$.
\end{theorem}
\begin{proof}
We have to show that the complex $\zz/n(q)$ is canonically
quasi-isomorphic in the etale topology to the sheaf ${\bf \mu}_n^{\oo
q}$. In the category of complexes of
sheaves with transfers of $\zz/n$-modules in the etale topology with
homotopy invariant cohomology sheaves $DM_{-}^{eff}(k,\zz/n,et)$ we have 
$$\zz/n(q)=(\zz/n(1))^{\oo{q}}$$
Since $\zz(1)={\bf G}_m[-1]$ in the etale topology we have
$\zz/n(1)=\mu_n$ which proves our claim.
\end{proof}
\begin{remark}\rm
A very detailed proof of this theorem which uses only the most basic
facts about motivic cohomology can be found in \cite{MW}.
\end{remark}

\begin{cor}
\llabel{c1follows} Assume that $H90(w,l)$ holds. Then for any $k$ of
characteristic not equal to $l$ and any $q\le w$ the norm residue map
\begin{equation}
\llabel{stat}
K_q^M(k)/l\sr H^{n}_{et}(Spec(k),{\bf \mu}_l^{\oo q})
\end{equation}
is an isomorphism.
\end{cor}
\begin{proof}
By Corollary \ref{svnew}(2) the homomorphisms
\begin{equation}
\llabel{fch}
H^{p,q}(X,\zz/l)\sr H^{p,q}_L(X,\zz/l)
\end{equation}
are isomorphisms for $p\le q\le w$. By Theorem \ref{sv0}, for $p=q$ and
$X=Spec(k)$ the homomorphism (\ref{fch}) is isomorphic to the
homomorphism (\ref{stat}).
\end{proof}
\begin{lemma}
\llabel{bkhilb} Assume that $H90({w},l)$ holds.  Then for any field
$k$, any $q\le {w}$ and any cyclic extension $E/k$ of degree $l$ the
sequence
$$K_q^M(E)\stackrel{1-\sigma}{\sr}
K_q^M(E)\stackrel{N_{E/k}}{\sr}K_q^M(k)$$
(where $\sigma$ is a generator of $Gal(E/k)$) is exact.
\end{lemma}
\begin{proof}
One verifies easily that this complex becomes exact after tensoring
with $\zz[1/l]$. It remains to show that it becomes eact after
tensoring with $\zz_{(l)}$.  Recall that for a smooth scheme $X$ over
$k$ we let $\zz_{tr}(X)$ denote the free sheaf with transfers
generated by $X$. Consider the complex of presheaves with transfers of
the form
$$0\sr \zz_{tr}(k)\sr \zz_{tr}(E)\stackrel{1-\sigma}{\sr}
\zz_{tr}(E)\sr \zz_{tr}(k)\sr 0$$
where the second arrow is the transfer map. Denote this complex with
the last $\zz_{tr}(k)$ placed in degree zero by $\uu{K}$. One can
easily see that it is exact in the etale topology and therefore 
$$Hom_{D}(\uu{K},{\bf
R}\pi_*(\pi^*(\zz(q)))[*])=Hom_{D_{et}}(\uu{K},\zz(q)[*])=0$$ 
where $D_{et}$ is the category of etale sheaves with transfers. Since
$$\uu{H}^{q+1}({\bf R}\pi_*(\pi^*(\zz(q))))=0,$$
the map
$$Hom_{D}(\uu{K},L(q)[q+2])\sr Hom_{D}(\uu{K},{\bf
R}\pi_*(\pi^*(\zz(q)))[q+2])$$
is injective and we conclude by Theorem \ref{sv} that
\begin{equation}
\llabel{iszero7}
Hom_{D}(\uu{K},\zz_{(l)}(q)[q+2])=0
\end{equation} 
We have
$$H^{p,q}(X,\zz)=Hom_{D}(\zz_{tr}(X),\zz(q)[p])$$
and in particular for a separable extension $F$ of $k$ we have
$$Hom_{DM}(\zz_{tr}(Spec(F)),\zz(q)[q])=\left\{
\begin{array}{ll}
K_q^M(F) &\mbox{\rm for $p=q$}\\
0 &\mbox{\rm for $p>q$}
\end{array}
\right.
$$
Our result follows now from (\ref{iszero7}) and the standard spectral
sequence which computes morphisms in a triangulated category from a
complex in terms of morphisms from its terms.
\end{proof}

\comment{The following result is the etale analog of Proposition \ref{indzar}.
\begin{proposition}
\llabel{indet} Let $k$ be a field and $k_0$ a subfield in $k$ such
that $k$ is separable over $k_0$. Let further $X/k$ be a smooth
scheme over $k$ and $X/k_0$ be the same scheme considered as a smooth
scheme over $k_0$. Then the natural morphism 
$$H^{p,q}_L(X/k_0,A)\sr H^{p,q}_L(X/k,A)$$
is an isomorphism.
\end{proposition}}

\subsection{Main theorem}
\label{fourth}\llabel{sec6}
In this section we prove the following theorem.
\begin{theorem}
\llabel{maintheorem} For any field $k$ and any $w\ge 0$ one has
$$H^{w+1,w}_L(Spec(k),\zz_{(2)})=0.$$
\end{theorem}
\begin{proof}
For $char(k)=2$ the theorem is proved in \cite[Theorem
8.6]{GL}. Assume that $char(k)\ne 2$. By induction on $w$ we may
assume that $H^{q+1,q}_L(k,\zz_{(2)})=0$ for all $q<w$.

Let $k$ be a field which has no extensions of degree prime to $2$ and
such that $K_{w}^M(k)$ is 2-divisible. By Lemma \ref{bkhilb} our
inductive assumption implies that $BK(w-1,2)$ holds. By Theorem
\ref{l4} we conclude that 
$$H^{w}_{et}(Spec(k),\zz/2)=0.$$
Together with Theorem \ref{sv0} this shows that the group
$H^{w+1,w}_{L}(k,\zz_{(2)})$ is torsion free.  On the other hand
$H^{w+1,w}_{L}(k,\qq)=0$ by Lemma \ref{ratLnew} and we conclude that
$H^{{w}+1,{w}}_L(Spec(k),\zz_{(2)})=0$.

For a finite extension $E/k$ of degree prime to $2$ the homomorphism
$$H^{{w}+1,{w}}_L(Spec(k),\zz_{(2)})\sr
H^{{w}+1,{w}}_L(Spec(E),\zz_{(2)})$$
is a monomorphism by the transfer argument. Thus in order to prove
$H90({w},2)$ it is sufficient to show that for any element $\uu{a}\in
K_{w}^M(k)$ there exists an extension $K_{\uu{a}}/k$ such that
$\uu{a}$ is divisible by $2$ in $K^M_{w}(K_{\uu{a}})$ and the
homomorphism
\begin{equation}
\llabel{sbm}
H^{{w}+1,{w}}_L(Spec(k),\zz_{(2)})\sr
H^{{w}+1,{w}}_L(Spec(K_{\uu{a}}),\zz_{(2)})
\end{equation}
is a monomorphism. Since any element in $K_{w}^M(k)$ is a sum of
symbols it is sufficient to construct $K_{\uu{a}}$ for $\uu{a}$ of the
form $(a_1,\dots,a_{w})$.

Let us show that the function field $K=k(Q_{\uu{a}})$ of the norm
quadric has the required properties. The fact that $\uu{a}$ becomes
divisible by $2$ in the Milnor K-theory of $K$ proved in Proposition
\ref{isspl0}. It remains to show that the map (\ref{sbm}) is
injective. We do it in two steps - first we prove that the kernel of
(\ref{sbm}) is covered by the group $H^{w+1,w}({\cal
X}_{\uu{a}},\zz_{(2)})$ and then that the later group is zero. For the
first step we will need the following three lemmas.
\begin{lemma}
\llabel{simpleneed}
Let $X$ be a non empty smooth scheme over $k$. Then the homomorphisms
\begin{equation}
\llabel{homo1}
H^{p,q}_L(Spec(k),\zz)\sr H^{p,q}_L(\check{C}(X),\zz)
\end{equation}
defined by the morphism $\check{C}(X)\sr Spec(k)$ are
isomorphisms for all $p,q\in\zz$. 
\end{lemma}
\begin{proof}
By definition, (see Appendix A) we can rewrite the homomorphism
(\ref{homo1}) as the homomorphism
$$Hom_{D_{et}}(\zz,\zz(q)[p])\sr
Hom_{D_{et}}(\zz(\check{C}(X)),\zz(q)[p])$$
where the morphisms are in the derived category of sheaves of abelian
groups in the etale topology on $Sm/k$. We claim that the morphism of
complexes $\zz(\check{C}(X))\sr \zz$ is a quasi-isomorphism. To check
this statement we have to verify that for any strictly henselian local
scheme $S$ the map of complexes of abelian groups
$\zz(\check{C}(X))(S)\sr \zz$ is a quasi-isomorphism. Since $X$ is
non-empty there exists a morphism $S\sr X$ and therefore, the
simplicial set $\check{C}(X)(S)$ is contractible (see the proof of
Lemma \ref{cc}). Therefore we have
$$\zz(\check{C}(X))(S)=\zz(\check{C}(X)(S))\cong \zz$$
\end{proof}
Let $K({w})$ be the complex of sheaves on $(Sm/k)_{Nis}$ defined in
Section \ref{first}. 
\begin{lemma}
\llabel{ctnew} 
Let $k$ be a field of characteristic not equal to $l$
and assume that $H90(w-1,l)$ holds. Then for any smooth scheme $X$ 
the map 
$${\bf H}^{*}(X,K(w)\oo\zz_{(l)})\sr
{\bf H^{*}}(X\times(\af-\{0\}),K(w)\oo\zz_{(l)})$$
defined by the projection $\af-\{0\}\sr Spec(k)$, is an isomorphism.
\end{lemma}
\begin{proof}
Recall from \cite{H2new} that for a functor $F$ from schemes to
abelian groups we denote by $F_{-1}$ the functor $X\mapsto
coker(F(X)\sr F(X\times{\bf A}^1-\{0\}))$ where the map is defined by
the projection. To prove the lemma we have to show that ${\bf
H}^{*}(X,K(w)\oo\zz_{(l)})_{-1}=0$. This is clearly equivalent to
checking that the standard map
$$H^{*,w}(X,\zz_{(l)})_{-1}={\bf H}^*(X,\zz_{(l)}(w))_{-1}\sr {\bf
H}^*_L(X,\zz_{(l)}(w))_{-1}$$
is an isomorphism. For a complex of sheaves with transfers $K$ there
is a complex $K_{-1}$ (defined up to a canonical quasi-isomorphism)
such that ${\bf H}^*(-,K)_{-1}={\bf H}^*(-,K_{-1})$. By
\cite[Proposition 4.34, p.124]{H2new}, if $K$ is a complex of sheaves
with transfers with homotopy invariant cohomology sheaves $\uu{H}^i$
then
$$\uu{H}^i(K_{-1})=(\uu{H}^i(K))_{-1}$$
Therefore, it is sufficient to check that the maps
$$\uu{H}^i(\zz_{(l)}(w))_{-1}\sr \uu{H}^i(L(w)\oo \zz_{(l)})_{-1}$$
are isomorphisms. Since both sides are zero for $i>w+1$ and since
$\uu{H}^i(K)$ are the sheaves associated with the presheaves $X\mapsto
{\bf H}^i(X,K)$ it remains to check that the maps
$${H}^{i,w}(X,\zz_{(l)})_{-1}\sr {\bf H}^i(X,L(w)\oo\zz_{(l)})_{-1}$$
are isomorphisms for $i\le w+1$. In this range the map 
$${\bf H}^i(X,L(w)\oo\zz_{(l)})\sr {\bf H}^i(X,{\bf
R}\pi_*(\pi^*(\zz_{l}(w))))=H^{i,w}_L(X,\zz_{(l)})$$
is an isomorphism and therefore it remains to check that the map
$$H^{i,w}(X,\zz_{(l)})_{-1}\sr H^{i,w}_L(X,\zz_{(l)})_{-1}$$ 
is an isomorphism for $i\le w+1$. Consider the commutative diagram
$$
\begin{CD}
H^{i-1,w-1}(X,\zz_{(l)}) @>>> H^{i-1,w-1}_L(X,\zz_{(l)})\\
@VVV @VVV\\
H^{i,w}(X,\zz_{(l)})_{-1} @>>> H^{i,w}_L(X,\zz_{(l)})_{-1} 
\end{CD}
$$
where the vertical arrows are defined by the multiplication with the
canonical class $\eta\in H^{1,1}(\af-\{0\},\zz)$. The upper horizontal
arrow is an isomorphism by our assumption that $H90(w-1,l)$ holds and
Corollary \ref{svnew}(1). The left hand side vertical arrow is an
isomorphism by the supension theorem \cite[Theorem 2.4]{Red}. It
remains to check that the right hand side vertical arrow is an
isomorphism. This we can verify separately for rational coefficients and
$\qq/\zz_{(l)}$-coefficients. In the former case the result follows
from Lemma \ref{ratLnew} and again \cite[Theorem 2.4]{Red}. In the
later case it follows from Theorem \ref{sv0} and the corresponding
result for the etale cohomology.
\end{proof}
\begin{lemma}
\llabel{nol} Assume that $H90(w-1,l)$ holds and let $k$ be a field of
characteristic not equal to $l$, $X$ a smooth scheme over $k$ and $U$
a dense open subscheme in $X$. Then the map
\begin{equation}
\llabel{homo9}
{\bf H}^{*}(X,K(w))\sr {\bf H}^{*}(U,K(w))
\end{equation}
is an isomorphism.
\end{lemma}
\begin{proof}
Considering $X$ to be a limit of smooth schemes
(possibly of greater dimension) over the subfield of constants in $k$
we may assume that $k$ is perfect. By obvious induction
it is sufficient to show that the statement of the lemma holds for
$U=X-Z$ where $Z$ is a smooth closed subscheme in $X$. Moreover, one
can easily see that it is sufficient to prove that for any point $z$
of $Z$ there exists a neighborhood $V$ of $z$ in $X$ such that 
$${\bf H}^{*}(V,K(w))\sr {\bf H}^{*}(V-V\cap Z,K(w))$$
is an isomorphism. This follows from the easiest case of the Gysin
distinguished triangle \cite[Proposition 3.5.4, p.221]{H3new} or
from the homotopy purity theorem \cite[Theorem 2.23, p.115]{MoVo} and
Lemma \ref{ctnew}.
\end{proof}
({\bf Proof of Theorem \ref{maintheorem} continues}) Let $u$ be an
element in the kernel of (\ref{sbm}). By Lemma \ref{simpleneed} it is
sufficient to show that the image $u'$ of $u$ in
$H^{{w}+1,{w}}_L({\cal X}_{\uu{a}}, \zz_{(2)})$ belongs to the image
of the group $H^{{w}+1,{w}}({\cal X}_{\uu{a}},\zz_{(2)})$ i.e. that
under our assumptions the image of $u'$ in the hypercohomology group
${\bf H}_{Nis}^{{w}+1}({\cal X}_{\uu{a}},K({w}))$ is zero.  Since $u$
becomes zero in the generic point of $Q_{\uu{a}}$ there exists a
nonempty open subscheme $U$ of $Q_{\uu{a}}$ such that the restriction
of $u'$ to $U$ is zero.  By Lemma \ref{nol} we conclude that the
restriction of $u'$ to $Q_{\uu{a}}$ is zero.

The canonical morphism $M(Q_{\uu{a}})\sr M({\cal X}_{\uu{a}})$ factors
through the morphism $M_{\uu{a}}\sr M({\cal X}_{\uu{a}})$ which is a
part of the distinguished triangle of Theorem \ref{Pmain}. Since
$M(Q_{\uu{a}})\sr M_{\uu{a}}$ has a section, our class $u'$ becomes
zero on $M_{\uu{a}}$ and by (\ref{maindt}) we conclude that it belongs
to the image of the group $Hom(M({\cal
X}_{\uu{a}})(2^{w-1}-1)[2^{w}-1],K(w)[w+1])$. Since $w>1$ this group
is zero by Lemma \ref{ctnew}.

It remains to show that $H^{w+1,w}({\cal X}_{\uu{a}},\zz_{(2)})=0$.
Consider the pointed simplicial scheme $\tilde{\cal
X}_{\uu{a}}=\wt{C}(Q_{\uu{a}})$ defined by the cofibration sequence
$$({\cal X}_{\uu{a}})_+\sr S^0\sr 
\tilde{\cal X}_{\uu{a}}\sr \Sigma^1_s(({\cal X}_{\uu{a}})_+)$$
Since $H^{p,q}(Spec(k),\zz)=0$ for $p>q$ the homomorphisms
$${H}^{p-1,q}({\cal X}_{\uu{a}},\zz)\sr
\wt{H}^{p,q}(\tilde{\cal X}_{\uu{a}},\zz)$$
defined by the third arrow of this sequence are isomorphisms for
$p-1>q$. Thus it is sufficient to verify that
$\tilde{H}^{{w}+2,{w}}(\tilde{\cal X}_{\uu{a}},\zz_{(2)})=0$. Since
there exists an extension of degree two $E/k$ such that
$Q_{\uu{a}}(Spec(E))\ne \emptyset$ Lemma \ref{van0a} implies that all
the motivic cohomology groups of $\tilde{\cal X}_{\uu{a}}$ have
exponent at most $2$. Thus it is sufficient to show that the image of
$\tilde{H}^{{w}+2,{w}}(\tilde{\cal X}_{\uu{a}},\zz_{(2)})$ in
$\tilde{H}^{{w}+2,{w}}(\tilde{\cal X}_{\uu{a}},\zz/2)$ is zero. Let
$u$ be an element of this image. Consider the composition of
cohomological operations $Q_{{w}-2}Q_{{w}-3}\dots Q_1$. It maps $u$ to
an element of
$$\tilde{H}^{2^{w},2^{{w}-1}}(\tilde{\cal
X}_{\uu{a}},\zz/2)={H}^{2^{w}-1,2^{{w}-1}}({\cal
X}_{\uu{a}},\zz/2)$$
\begin{lemma}
\llabel{inttoint}
Let $a$ be a class in $H^{*,*}({\cal X},\zz/l)$ which the image of an
integral class. Then $Q_i(a)$ is the image of an integral
class for any $i$.
\end{lemma}
\begin{proof}
By \cite[Proposition 13.6]{Red} we have $Q_i=\beta q_i\pm
q_i\beta$. Since $a$ is the image of an integral class we have
$\beta(a)=0$ and $Q_i(a)=\beta q_i(a)$. On the other hand, the
Bockstein homomorphism $\beta$ can be written as the composition
$$\wt{H}^{*,*}(-,\zz/l)\stackrel{B}{\sr}\wt{H}^{*+1,*}(-,\zz)\sr
\wt{H}^{*+1,*}(-,\zz/l)$$
where the first map is the connecting homomorphism for the exact
sequence $0\sr \zz\sr \zz\sr\zz/l\sr 0$ and the second map is the
reduction modulo $l$. Therefore, any element of the form $\beta(x)$ is
the image of an integral class.
\end{proof}
Lemma \ref{inttoint} implies that this element belongs to the image of
the corresponding integral motivic cohomology group. Therefore by
Theorem \ref{rostvan} it is zero. It remains to verify that the
composition $Q_{{w}-2}Q_{{w}-3}\dots Q_1$ is a monomorphism i.e that
the operation $Q_i$ acts monomorphically on the group
$\tilde{H}^{{w}-i+2^{i}-1,{w}-i+2^{i-1}-1}(\tilde{\cal
X}_{\uu{a}},\zz/2)$ for $i=1,\dots,{w}-2$. By Corollary \ref{exitcor},
the motivic Margolis homology of $\wt{\cal X}_{\uu{a}}$ are zero for
$Q_i$ with $i\le {w}-1$. Therefore, the kernel of $Q_i$ on this group
is covered by the image of $\tilde{H}^{{w}-i,{w}-i}(\tilde{\cal
X}_{\uu{a}},\zz/2)$. The later group is zero by the inductive
assumption that $H90({w}-1,2)$ holds, Lemma \ref{simpleneed} and
Corollary \ref{sv}(2).
\end{proof}
Theorem \ref{maintheorem} is proved.
\begin{cor}
\llabel{MC}
Let $k$ be a field of characteristic not equal to $2$. Then the norm residue
homomorphisms $K_{w}^M(k)/2\sr H^{w}_{et}(k,\zz/2)$ are isomorphisms for all
${w}\ge 0$.
\end{cor}
\begin{cor}
For any field $k$ of characteristic not equal to $2$ and any $q\ge 0$
the etale cohomology group $H^{q+1}_{et}(k,\zz_2(q))$ is torsion free.
\end{cor}
The following nice corollary of Theorem \ref{maintheorem} is due to S. Bloch.
\begin{cor}
Let $\alpha\in H^i(X,\zz)$ be a
2-torsion element in the integral cohomology of a complex algebraic
variety $X$. Then there exists a divisor $Z$ on $X$ such that
the restriction of $\alpha$ to $X-Z$ is zero.
\end{cor}
\begin{proof}
Since $2\alpha=0$, $\alpha$ is the image of an element $\alpha'$ in
$H^{i-1}(X,\zz/2)$ with respect to the 
Bockstein homomorphism $\beta:H^{i-1}(X,\zz/2)\sr H^i(X,\zz)$. By Corollary
\ref{MC} there exists a dense open subset $U=X-Z$ of $X$ such that
$\alpha'$ on $U$ is in the image of the canonical map $({\cal
O}^*(U))^{\oo (i-1)}\sr  H^{i-1}(U,\zz/2)$. Since this map factors
through the integral cohomology group we have $\beta(\alpha')=0$ on
$U$.
\end{proof}
Denote by $B/n(q)$ the canonical truncation at the cohomological level
$q$ of the complex of sheaves ${\bf R}\pi_*(\pi^*({\bf \mu}_n^{\oo
q}))$ where $\pi:(Sm/k)_{et}\sr (Sm/k)_{Nis}$ is the usual
morphism of sites. 
\begin{theorem}
\llabel{alaB} For a smooth variety $X$ and any $n>0$ there are
canonical isomorphisms
$$H^{p,q}(X,\zz/2^n)={\bf H}_{Nis}^p(X,B/2^n(q))$$
In particular, for any $X$ as above $H^{p,q}(X,\zz/2^n)=0$ for $p<0$.
\end{theorem}

In \cite{BL} Spencer Bloch and Stephen Lichtenbaum constructed the
{\em motivic spectral sequence} which starts from Higher Chow groups
of a field and converges to its algebraic K-theory. Their construction
was reformulated in much more natural terms and extended to all smooth
varieties over fields in \cite{FS}. Combining the motivic spectral
sequence with Theorem \ref{alaB} and using \cite{comparison} to
identify motivic cohomology of \cite{FS} with the motivic cohomology
of this paper we obtain the following result.
\begin{theorem}
\llabel{alaQL} Let $k$ be a field of characteristic $\ne 2$. Then for
all smooth $X$ over $k$ and all $n>0$ there exists a natural spectral
sequence with the $E_2$-term of the form
$$E_2^{p,q}=\left\{
\begin{array}{ll}
H^{p-q}_{et}(X,\zz/2^n(q))&\mbox{\rm for $p,q\le 0$}\\
0&\mbox{\rm otherwise}
\end{array}
\right.
$$
which converges to $K_{-p-q}(X,\zz/2^n)$ (K-theory with
$\zz/2^n$-coefficients). 
\end{theorem}
\begin{cor}
Let $X$ be a complex algebraic variety of dimension $d$. Then the
canonical homomorphisms $K_i^{alg}(X,\zz/2^n)\sr K_i(X({\bf
C}),\zz/2^n)$ are isomorphisms for $i\ge d-1$ and monomorphisms for
$i=d-2$.
\end{cor}

\subsection{Appendix A. Hypercohomology of pointed simplicial sheaves}
\llabel{sec0}
We recall here some basic notions related to the hypercohomology of
simplicial sheaves.  Let $T$ be a site with final object $pt$. For the
purpose of the present paper one may assume that $T$ is the category
of smooth schemes over a field $k$ with the Nisnevich or the etale
topology and $pt=Spec(k)$. Denote by $AbShv(T)$ the category of
sheaves of abelian groups on $T$. For an object $X$ of $T$ let
$\zz(X)$ be the sheaf characterized by the property that
\begin{equation}
\llabel{adj}
Hom(\zz(X),F)=F(X)
\end{equation}
for any sheaf of abelian groups $F$. Note that $\zz(Spec(k))$ is the
constant sheaf $\zz$.

For a pointed simplicial object ${\cal X}$, $\zz({\cal X})$ is a
simplical sheaf of abelian groups and we may form a complex $\zz({\cal
X})_*$ by taking the alternating sums of boundary maps. For a complex
of sheaves of abelian groups $K$ we define the hypercohomology of
$\cal X$ with coefficients in $K$ by the formula
$${\bf H}^n({\cal X},K):=Hom_{D}(\zz({\cal X}),K)$$
where $D=D(AbShv(T))$ is the derived category of complexes of sheaves
of abelian groups on $T$. The formula (\ref{adj}) together with the
standard method of computing cohomology by means of injective
resolutions implies that for ${\cal X}=X$ an object of $T$ our
definition agrees with the usual one.

Let now $Y$ be a sheaf of sets on $T$. Then we can define $\zz(Y)$ as
the free sheaf of abelian groups generated by $Y$ such that for every
sheaf of abelian groups $F$ one has $Hom(\zz(Y),F)=Hom(Y,F)$ where on
the left hand side we have morphisms of sheaves of abelian groups and
on the right hand side morphisms of sheaves of sets. The Yoneda Lemma
shows that for an object $X$, $\zz(X)=\zz(h_X)$ where $h_X$ is the
sheaf of sets represented by $X$. Therefore, we may immediately extend
our definition of hypercohomology groups to simplicial sheaves in a
way which agrees with the definition for simplicial objects of $T$ on
representable sheaves. Starting from this point we consider objects of
$T$ to be a particular type of sheaves of sets.

If ${\cal X}$ is a pointed simplicial sheaf the distinguished point
defines a homomorphism ${\bf H}^*({\cal X},K)\sr {\bf H}^*(pt,K)$
which has a canonical section defined by the projection ${\cal X}\sr
pt$. We define reduced hypercohomology of ${\cal X}$ by the
formula:
$$\wt{H}^*({\cal X},K)=ker({\bf H}^*({\cal X},K)\sr {\bf
H}^*(pt,K))$$
Alternatively, we can define $\wt{\zz}({\cal X})$ setting
$$\wt{\zz}({\cal X}):=coker(\zz\sr{\zz}({\cal X}))$$
where the morphism is defined by the distinguished point. Then
$$\wt{H}^n({\cal X},K)=Hom_D(\wt{\zz}(X),K[n])$$
If ${\cal X}_+$ denotes the simplicial sheaf ${\cal X}\coprod pt$
pointed by the canonical embedding $pt\sr {\cal X}\coprod
Spec(k)$, then
$$H^{*,*}({\mcal X},-)=\tilde{H}^{*,*}({\mcal X}_+,-)$$
There is a standard spectral sequence with the $E_1$-term consisting
of the hypercohomology of the terms ${\cal X}_i$ of $\cal X$ with
coefficients in $K$ which tries to converge to hypercohomology of
$\cal X$. However, this spectral sequence is of limited use since
$\zz({\cal X})_*$ is unbounded on the left and the convergence
properties of this spectral sequence are uncertain. The proof of the
following proposition contains the trick which allows one get around
this problem.
\begin{proposition}
\llabel{vangen} Let $n\ge 0$ be an integer and ${\cal X}$ and $K$ be
such that 
$${\bf H}^m({\cal X}_j,K)=0$$
for all $j\ge 0$ and $m\le n$. Then ${\bf H}^m({\cal X},K)=0$ for
$m\le n$.
\end{proposition}
\begin{proof}
Let $s_{\le i}\zz({\cal X})_{*}$ be the subcomplex of $\zz({\cal
X})_{*}$ such that 
$$
s_{\le i}\zz({\cal X})_{j}=\left\{
\begin{array}{ll}
\zz({\cal X})_{j}&\mbox{\rm for $j\le i$}\\
0&\mbox{\rm for $j > i$}
\end{array}
\right.
$$
Since the complexes $s_{\le i}\zz({\cal X})_{*}$ are bounded the usual
argument shows that $Hom_{D}(s_{\le i}\zz({\cal X})_{*}, K[m])=0$ for
$m\le n$ and all $i$. On the other hand one has an exact sequence of
complexes
$$0\sr \oplus_i s_{\le i}\zz({\cal X})_{*}\sr \oplus_i s_{\le
i}\zz({\cal X})_{*} \sr \zz({\cal X})_{*}\sr 0$$
which expresses the fact that $\zz({\cal X})_{*}=colim_is_{\le
i}\zz({\cal X})_{*}$. This short exact sequence defines a long exact
sequence of groups of morphisms in the derived category which shows
that ${\bf H}^m({\cal X},K)=0$ for $m\le n$.
\end{proof}
For a morphism of simplicial sheaves $f:{\cal X}\sr {\cal X}'$ let
$cone(f)$ denote the simplicial sheaf such that for a smooth scheme
$U$ one has 
$$cone(f):U\mapsto cone({\cal X}(U)\sr {\cal X}'(U))$$
\begin{lemma}
\llabel{longex}
The sequence
$${\cal X}\stackrel{f}{\sr} {\cal X}'\sr cone(f)$$
defines a long exact sequence of hypercohomology groups of the form
$$\dots\sr \wt{H}^{n}(cone(f),K)\sr \wt{H}^{n}({\cal X},K)\sr
\wt{H}^{n}({\cal X}',K)\sr$$
$$\sr\wt{H}^{n+1}(cone(f),K)\sr\dots$$
\end{lemma}

\subsection{Appendix B. \v Cech simplicial schemes}
We define here for any smooth scheme $X$ a simplicial scheme
$\check{C}(X)$ such that the map $\pi:\check{C}(X)\sr Spec(k)$ is a
weak equivalence if and only if $X$ has a rational point.  This map
defines an isomorphism in motivic cohomology with
$\zz_{(n)}$-coefficients if and only if $X$ has a 0-cycle of degree
prime to $n$. In particular, the motivic cohomology of
$\wt{C}(X):=cone(\pi)$ with $\zz/n$-coefficients provide obstructions
to the existence of 0-cycles on $X$ of degree prime to $n$.
\begin{definition}
For a variety $X$ over $k$, $\check{C}(X)$ is the
simplicial scheme with terms $\check{C}(X)_n=X^{n+1}$ and face and
degeneracy morphisms given by partial projections and diagonals
respectively.
\end{definition}
\begin{lemma}
\llabel{cc}
Let $X,Y$ be smooth schemes over $k$ such that 
$$Hom(X,Y)\ne \emptyset$$
Then the projection $\check{C}(Y)\times X\sr X$
is a simplicial weak equivalence. 
\end{lemma}
\begin{proof}
Let $U$ be a smooth scheme over $k$. Then the simplicial set
$\check{C}(Y)(U)$ is the simplex generated by the set $Hom(U,Y)$ which
is contractible if and only if $Hom(U,Y)\ne\emptyset$. This implies
immediately that for any $U$ the map of simplicial sets
$$(\check{C}(Y)\times X)(U)\sr X(U)=Hom(X,U)$$ is a weak
equivalence. 
\end{proof}
Let $\wt{C}(X)$ denote the unreduced suspension of $\check{C}(X)$
i.e. the cone of of the morphism $\check{C}(X)_+\sr Spec(k)_+$. 
\begin{lemma}
\llabel{van0a} Let $Y$ be a smooth scheme which has a rational point
over an extension of $k$ of degree $n$. Then
$n\wt{H}^{*,*}(\wt{C}(Y),\zz)=0$.
\end{lemma}
\begin{proof}
Let $E$ be an extension of $k$ and $Y_E=Y\times_{Spec(k)}Spec(E)$
considered as a smooth scheme over $E$. Then there are homomorphisms
$$\wt{H}^{*,*}(\wt{C}(Y),\zz)\sr \wt{H}^{*,*}(\wt{C}(Y_E),\zz)$$
and
$$\wt{H}^{*,*}(\wt{C}(Y_E),\zz)\sr \wt{H}^{*,*}(\wt{C}(Y),\zz)$$
and the composition of the first one with the second is multiplication
by $deg(E/k)$. Let $E$ be an extension of degree $n$ such that
$Y(E)\ne \emptyset$. Since $Y(E)\ne \emptyset$ the pointed sheaf
$\wt{C}(Y_E)$ is contractible and therefore,
$$\wt{H}^{*,*}(\wt{C}(Y_E),\zz)=0.$$
We conclude that $n\wt{H}^{*,*}(\wt{C}(Y),\zz)$.
\end{proof}

\end{document}